%% file: UfMJ.tex
\documentclass[review,hidelinks,onefignum,onetabnum]{siamart220329}

\input{ex_shared}

\ifpdf
\hypersetup{
  pdftitle={Unification of Lagrangian staggered-grid hydrodynamics and cell-centered hydrodynamics in one dimension},
  pdfauthor={Xihua Xu}
}
\fi

\begin{document}

\maketitle

\begin{abstract}
This paper focuses on the novel scheme to unify both Lagrangian staggered-grid and cell-centered hydrodynamic methods in one dimension. The scheme neither contains empirical parameters nor solves the Riemann problem. It includes two key points: one is the relationship between pressure and velocity, and the other is Newton's second law. The two methods that make use of this scheme satisfy the entropy condition and are conservative in total mass, momentum, and energy. Numerical results show the robustness and accuracy of both methods.  
\end{abstract}



\section{Introduction}
For a number of years, both Lagrangian staggered-grid hydrodynamic (SGH) and cell-centered hydrodynamic (CCH) methods have received much attention due to their natural properties. One of these properties is that the moving reference frame is advantageous for computing contact discontinuities, shock waves, and material interfaces. That is the reason why the two methods are widely applied in the simulation of multi-material and free surface flows, and have a broad range of applications in fields such as steam explosion, inertial confinement fusion, and high energy density physics.

The SGH method employs a staggered discretization in which position, velocity, and kinetic energy are centered at points while density, pressure, and specific internal energy are defined within cells. This method originated with Von Neumann and Richtmyer \cite{von_neumann_method_1950}, and Wilkins \cite{wilkins_calculation_1963} extended it to multi-dimensions. Many researchers have improved the accuracy and the robustness of the SGH method, such as compatible staggered discretization \cite{caramana_construction_1998}, the anti-hourglass correction \cite{caramana_elimination_1998}, the edge viscosity \cite{caramana_formulations_1998}, the tensor viscosity \cite{campbell_tensor_2001} and their references.

In the CCH method, all conserved quantities (mass, momentum, and energy) are defined at the cell-center. Based on the Godunov \cite{godunov_1959} method, many developments have been made in order to calculate the nodal velocity and interface fluxes coherently, such as the least-squares procedure \cite{dukowicz_vorticity_1992}\cite{cline_caveat-gt:_1990}, the average of the two tangential parts \cite{cheng_high_2007}, the two-dimensional nodal solvers \cite{despres_lagrangian_2005}\cite{maire_cell-centered_2007}\cite{burton_cell-centered_2013}, and three-dimensional nodal solvers \cite{carre_cell-centered_2009}\cite{georges_3d_2016}. The latter multi-dimensional nodal solvers satisfy the well-known geometric conservation law (GCL) \cite{thomas_geometric_1978} and construct a consistent way to determine the vertex velocity and the numerical flux at the interface. The readers interested in both SGH and CCH numerical methods might find a more detailed presentation in \cite{benson_computational_1992}\cite{loubere_chapter_2016}\cite{margolin_artificial_2022}. However, the intrinsic relationship between the SGH and CCH methods has to be discovered in depth. 

Under the viewpoint of shock-capturing schemes, one way to build the relationship between the SGH and CCH methods is to link the artificial viscosity \cite{wilkins_use_1980} with the Riemann solver \cite{dukowicz_general_1985}. In \cite{christensen_godunov_1990}, Christensen found that the artificial viscosity is equivalent to the HLL approximate Riemann solver under certain assumptions, which implied there was a potential synergy between the two numerical methods. Luttwak and Falcovitz \cite{luttwak_staggered_2005} applied an uniaxial tensor pseudo-viscosity to provide necessary dissipation at shocks to build the Staggered Mesh Godunov scheme. Burbeau-Augoula \cite{burbeau-augoula_node-centered_2010} introduced an extra degree of freedom, the fluid velocity within cells, to establish a link between the CCH and SGH methods, motivated by recent progress in CCH schemes \cite{despres_lagrangian_2005}\cite{maire_cell-centered_2007},  This extra degree of freedom is coupled to the nodal velocity by constructing a linear velocity vector field approximation with frame-invariant limitation. Later, a series of compatible Lagrangian discretizations and associated Riemann solver-based artificial viscosities were developed for the SGH method. Maire et al. \cite{maire_staggered_2011} employed the concept of sub-cell discretization, and the artificial viscosity force was derived by invoking Galilean invariance and thermodynamic consistency. Then Loubère et al. \cite{loubere_3d_2013} extended this approach to obtain a 3D frame-invariant vector limitation. Morgan et al. \cite{morgan_lagrangian_2014} incorporated the multidirectional Riemann-like problem into SGH since this Riemann-like solution was very robust against mesh instabilities \cite{burton_cell-centered_2013}. Xu et al. \cite{xu_godunov-type_2020} derived a general form of the artificial viscosity, which was linked with the Generalized Riemann Invariant relation to determine the viscosity coefficients. Although these methods all employed the Riemann solver to approximate the artificial viscosity on a staggered grid and proved to be very robust and accurate, little attention has been paid to the unification of SGH and CCH methods without increasing degrees of freedom.

The primary aim of this paper is to develop a unified framework for Lagrangian SGH and CCH methods in one dimension. The unified framework applies the same scheme to discretize the corresponding set of gas dynamic equations on a staggered-grid and cell-centered mesh, intended initially for the same purposes. The main benefit of the new scheme is the elimination of empirical parameters, which are often contained in the SGH method. Moreover, this scheme does not solve the Riemann problem that is frequently presented in the CCH method. The scheme includes two key points: one is the relationship between pressure and velocity, and the other is Newton's Second Law. The first key point is based on the fact that both of the Hugoniot and isentropic curves have the same first- and second-order derivatives from a state \cite{richard_supersonic_1976}. Using the Taylor expansion, the relationship between pressure and specific volume obtains the same expression whether it is obtained through a shock or an isentropic process. This can be viewed as a minor jump from one state to another. The characteristic time of a cell is introduced as the time for the sound speed to travel through the target cell. According to the motion of vertices in the characteristic time, the link between specific volume and velocity is also constructed. The two relationships provides the first key point. A new finding that the nodal acceleration is zero in each time step is employed for the CCH method. Instead of an exact/approximate Riemann solver, Newton’s second law can be applied to calculate the velocity and pressure. Once the velocity and pressure are known, then gas dynamic equations can be updated. Thus, this scheme can be interpreted as a unification of Lagrangian SGH and CCH methods. Both SGH and CCH methods, using the same scheme, maintain the conservation of total mass, monument, and energy and satisfy the entropy conditions.

This paper is organized as follows. The second section introduces the system of one-dimensional hydrodynamic equations, notions, and the characteristic time of a target cell and reviews the properties of Hugonius and isentropic curves. In the third section, the discretization of both SGH and CCH methods is established by using the same scheme, which contains the three steps. The main properties, such as conservatively and compatibility with the entropy condition, are shown in details. Finally, numerical results are provided to verify the accuracy and robustness of the two methods.

\section{Basic equations}
In a classical Lagrangian framework, the system of hydrodynamic equations describing the motion of a compressible gas is given by
\begin{eqnarray}
&&\frac{\partial\tau}{\partial t} - \frac{\partial u}{\partial x} = 0, \label{mass}\\
&&\frac{\partial{u}}{\partial t} + \frac{\partial P}{\partial x} = 0,  \label{mum}\\
&&\frac{\partial{\varepsilon} }{\partial t} + P\frac{\partial u}{\partial x} = 0, \label{slsenin} ~ \quad\hbox{(for SGH)} \label{enineqn} \\
&&\frac{\partial{E} }{\partial t} + \frac{\partial Pu}{\partial x} = 0, \label{clseng}\quad\hbox{(for CCH)}\label{engeqn} 
\end{eqnarray}
where $\tau$ stands for the specific volume, $\rho = \frac{1}{\tau}$ the density, $P$ the pressure, $u$ the velocity, and $E$ the energy, with $E = \varepsilon + \frac{1}{2}u^2$ that is the sum of specific internal energy $\varepsilon$ and kinetic energy $\frac{1}{2}u^2$. The previous system is equipped with a thermodynamic closure, Equation of State (EOS)
\begin{equation}
\label{eos}
P=P(\rho,\varepsilon).
\end{equation} 

The Lagrangian motion of a point $x$ is described by the trajectory equations
\begin{eqnarray}
\frac{dx(t)}{dt}= u(x(t),t),\quad x(0) = x_0. \label{dxu}
\end{eqnarray}
where ${x}_0$ is the initial position of the point. This equation expresses the Lagrangian motion of any point initially located at position $x_0$.

In the SGH method, one may discretize the equations (\ref{mass}), (\ref{mum}) and (\ref{enineqn}). Velocity and kinetic energy are centered at points, while other variables (density and internal energy) are defined at the cell-center. In the CCH method, the equations (\ref{mass}), (\ref{mum}) and (\ref{engeqn}) will be discretized. All conserved quantities, including cell velocity, are cell-centered. The nodal velocity is undetermined.

\subsection{Notation and Definition}
\begin{figure}[!h]
\centering
	\includegraphics[scale=1.2]{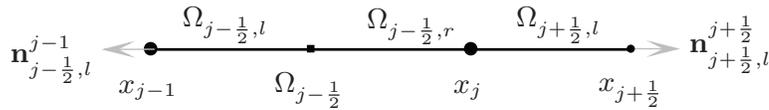}
	\caption{Notation for one-dimensional scheme description. 
	 \label{mesh1d}}	
\end{figure}

The computational domain $\Omega$ is a line segment of length $L$ and divided into $N$ cells denoted as $\Omega_{j-\frac{1}{2}},(j=1,...,N)$.  As shown in Figure \ref{mesh1d}, the nodes are labeled as $x_j,(j=1,...,N+1)$. The midpoints of each cell are defined by $x_{j-\frac{1}{2}},(j=1,...,N)$. Using the midpoints, each cell can be split into two sub-cells, denoted as  $\Omega_{j-\frac{1}{2},l}$ and $\Omega_{j-\frac{1}{2},r}$ where the subscript $l(r)$ means left(right) side sub-cell with respect to the cell.  

The notations allow for the definition of the unit outward normals of each cell and sub-cell. The unit outward normal of the cell $\Omega_{j-\frac{1}{2}}$ at the point $x_{j-1}$ is defined as $\textbf{n}^{j-1}_{j-\frac{1}{2}}$, the sub-cell $\Omega_{j+\frac{1}{2},l}$ at the point $x_{j+\frac{1}{2}}$ is $\textbf{n}^{j+\frac{1}{2}}_{j+\frac{1}{2},l}$, the sub-cell $\Omega_{j-\frac{1}{2},l}$ at the point $x_{j-1}$ is $\textbf{n}^{j-1}_{j-\frac{1}{2},l}$. Note that the subscripts denote the index of cells/sub-cells, and the superscripts are the index of nodes. The following identities shall be useful and hold true  
\begin{eqnarray}
\label{np1}
\textbf{n}^{j}_{j-\frac{1}{2}} = -\textbf{n}^{j}_{j+\frac{1}{2}}, \quad
\textbf{n}^{j-\frac{1}{2}}_{j-\frac{1}{2},l} = -\textbf{n}^{j-\frac{1}{2}}_{j-\frac{1}{2},r}, \quad 
\textbf{n}^{j}_{j-\frac{1}{2},r} = \textbf{n}^{j}_{j-\frac{1}{2}}.
\end{eqnarray}

For the sake of convenience, assuming that index growth is positive yields
\begin{eqnarray}
\label{np11}
\textbf{n}^{j+\frac{1}{2}}_{j+\frac{1}{2},l} = 1,\quad \textbf{n}^{j-1}_{j-\frac{1}{2}} = -1.
\end{eqnarray}

Moreover, remarking that the volumes of cell $\Omega_{j+\frac{1}{2}}$ and sub-cell $\Omega_{j+\frac{1}{2},l}(\Omega_{j+\frac{1}{2},r})$ can be expressed as 
\begin{eqnarray}
V_{j+\frac{1}{2}} = x_{j+1}-x_{j},\ V_{j+\frac{1}{2},l} = x_{j+\frac{1}{2}}-x_{j},\
V_{j+\frac{1}{2},r} = x_{j+1}-x_{j+\frac{1}{2}}.
\end{eqnarray}
By summation of Lagrangian sub-cell masses, one defines the masses of cells and nodes as 
\begin{eqnarray}
m_{j+\frac{1}{2}} &=& \rho_{j+\frac{1}{2}} V_{j+\frac{1}{2}}=\rho_{j+\frac{1}{2}}\left(x_{j+1}-x_j\right) \nonumber \\
m_{j}  &=&  m_{j-\frac{1}{2},r} + m_{j+\frac{1}{2},l} = \rho_{j-\frac{1}{2}} V_{j-\frac{1}{2},r} + \rho_{j+\frac{1}{2}} V_{j+\frac{1}{2},l} \nonumber \\
 &=& \rho_{j-\frac{1}{2}} \left(x_j - x_{j-\frac{1}{2}}\right) + \rho_{j+\frac{1}{2}}\left(x_{j+\frac{1}{2}}-x_j\right)\nonumber
\end{eqnarray}
where $m_{j-\frac{1}{2},r}$ and $m_{j+\frac{1}{2},l}$ are the masses of sub-cells $\Omega_{j-\frac{1}{2},r}$ and $\Omega_{j+\frac{1}{2},l}$, respectively. Since the cells move together with the fluid in the Lagrangian method, the masses of cells and nodes are constants as initial constants.

At the end of this subsection, the \underline{\textbf{characteristic time}} of a cell is introduced as the time for the sound speed to travel through the cell
\begin{eqnarray}
\delta{t} = \frac{l}{c}, \label{charatime}
\end{eqnarray}
where $l$ is the characteristic length of the cell and $c$ is the local sound speed. Each cell has its own characteristic time on account of its different geometric scales and local sound speeds. 

The characteristic time is one of the essential properties of a cell. This means that a small perturbation has traveled across the whole cell at the speed of sound, driving changes in certain physical quantities such as density, pressure, volume, and so on. This new finding will contribute to the understanding of the relationship between specific volume and velocity.

\subsection{Hugonius and isentropic curves}
In this section, some preliminary facts, helpful for the exposition, will be reviewed. The main properties of Hugonius and isentropic curves will be revisited from the well-known book \cite{richard_supersonic_1976}.

If the state $(\tau, P)$ reached from a certain initial state $(\tau_0,P_0)$ is compressed by a shock, the two states satisfy the following Hugoniot relation 
\begin{eqnarray}
H(\tau,P) = \varepsilon(\tau,P) - \varepsilon(\tau_0,P_0) +
\frac{1}{2}(\tau-\tau_0)(P-P_0)=0. \label{HugoniotRelation}
\end{eqnarray}
The pressure $P$ is uniquely defined by the specific volume $\tau$. This function can be denoted as $P=P_H(\tau)$ and yields $P_0=P_H(\tau_0)$. 

In the $(\tau, P)$ plane, the graph of the Hugoniot relation is a curve passing through the point $(\tau_0, P_0)$. This curve is called the Hugoniot curve. Along the Hugoniot curve, all the thermodynamic parameters are functions of $\tau$. Differentiating the (\ref{HugoniotRelation}) along the Hugoniot curve
\begin{eqnarray}
\label{dh}
\frac{d\varepsilon_H}{d\tau} + \frac{1}{2}(P_H+P_0) +
\frac{1}{2}(\tau-\tau_0)\frac{dP_H}{d\tau} =0.
\end{eqnarray}
By using the Gibbs formula, $TdS=d\varepsilon+Pd\tau$, the formula (\ref{dh}) becomes
\begin{eqnarray}
T_H\frac{dS_H}{d\tau} - \frac{1}{2}(P_H-P_0) +
\frac{1}{2}(\tau-\tau_0)\frac{dP_H}{d\tau} =0. \nonumber
\end{eqnarray}
Setting $\tau=\tau_0$ in the above equation reads
\begin{eqnarray}
T_H(\tau_0)\frac{dS_H(\tau_0)}{d\tau}=0 \quad\hbox{or}\quad
\frac{dS_H(\tau_0)}{d\tau}=0. \nonumber
\end{eqnarray}

Since $S_H(\tau)=S(\tau, P_H(\tau))$ holds true along the Hugoniot curve. Differentiating this equation with responding to $\tau$ and then setting $\tau=\tau_0$ leads to
\begin{eqnarray}
0=\frac{dS_H(\tau_0)}{d\tau}=\frac{\partial
S(\tau_0,P_0)}{\partial \tau}+\frac{\partial
S(\tau_0,P_0)}{\partial P}\frac{dP_H(\tau_0)}{d\tau}. \nonumber
\end{eqnarray}
Therefore, the following identity is derived
\begin{eqnarray}
\frac{dP_H(\tau_0)}{d\tau}=-\frac{\partial
S(\tau_0,P_0)}{\partial \tau}/\frac{\partial
S(\tau_0,P_0)}{\partial P}. \label{dph_dtau0}
\end{eqnarray}

Since $P$ can also be considered as a function of $\tau$ and $S$, i.e., $P=P(\tau, S)$. Substituting this equation into the equation $S=S(\tau, P)$ yields $S=S(\tau, P(\tau, S))$. Differentiating this equation with responding to $\tau$ gives
\begin{eqnarray}
0=\frac{\partial S}{\partial \tau}+\frac{\partial S}{\partial P}
\frac{\partial P}{\partial \tau}, \nonumber
\end{eqnarray}
this allows to deduce that
\begin{eqnarray}
\frac{\partial S}{\partial \tau}/\frac{\partial S}{\partial P}=
- \frac{\partial P}{\partial \tau} =\rho^2 c^2, \nonumber
\end{eqnarray}
Consequently, it is straightforward to check that (\ref{dph_dtau0}) can be rewritten as
\begin{eqnarray}
\label{ptau1}
\frac{dP_H(\tau_0)}{d\tau}=-\rho_0^2 c_0^2.
\end{eqnarray}

The above equation (\ref{ptau1}) can be shown as a linear approximation along the Hugoniot curve from the state $(\tau_0,P_0)$. To get a quadratic approximation, in the case of an ideal gas EOS, the second-order derivative reads
\begin{equation}
\label{ptau2}
\frac{d^2P_{H}(\tau_0)}{d{\tau}^2} = (\gamma+1)\rho_0^3c_0^2,
\end{equation}
where $\gamma$ is the adiabatic index. It is important to emphasize here that the equation (\ref{ptau2}) is only valid for the ideal gas EOS. For other equations of states, the second-order derivatives have different expressions, see \cite{wilkins_calculation_1963} \cite{richard_supersonic_1976}.

Since both Hugoniot and isentropic curves have the same first- and second-order derivatives from a state $(\tau_0,P_0)$, details in \cite{richard_supersonic_1976}, the equations (\ref{ptau1}) and (\ref{ptau2}) can be recast as
\begin{eqnarray}
\label{ptau1ptau2}
\frac{dP_{H}(\tau_0)}{d{\tau}}  =\frac{dP_{S}(\tau_0)}{d{\tau}}=-\rho_0^2 c_0^2, \quad 
\frac{d^2P_{H}(\tau_0)}{d{\tau}^2}  =\frac{d^2P_{S}(\tau_0)}{d{\tau}^2}(\gamma+1)\rho_0^3c_0^2 
\end{eqnarray}

The above equations (\ref{ptau1ptau2}) show that the state $(\tau,P)$ reached from an initial state $(\tau_0,P_0)$ by a shock wave or an isentropic process can be approximated by the following Taylor expansion
\begin{eqnarray}
\label{pptautau}
P &=& P_0 + \frac{dP(\tau_0)}{d\tau}(\tau-\tau_0) + \frac{1}{2}\frac{d^2P(\tau_0)}{d\tau^2}(\tau-\tau_0)^2   \nonumber \\[4mm]
&=& P_0 - \rho_0^2a^2_0 (\tau-\tau_0) +  \frac{\gamma+1}{2}\rho_0^3c_0^2(\tau-\tau_0)^2,
\end{eqnarray}
where a quadratic approximation is taken.

The relationship between pressure and specific volume has been built. This relationship implies that a minor jump from one state to another, caused by a shock wave or an isentropic process, can be calculated using the same expression, see (\ref{pptautau}). 

The formula (\ref{pptautau}) is similar to the so-called artificial viscosity \cite{caramana_formulations_1998}, where the force is written as a pressure contribution plus a tensorial viscous contribution. A new form of artificial viscosity is derived using the equivalence of lower-order derivatives along both Hugoniot and isentropic curves. 

\section{Discretization of the equations}
In this section, the same scheme is used to discretize the gas dynamics system in both SGH and CCH methods. In general, the dissipation of kinetic energy into internal energy through shock waves is ensured by an artificial viscosity term in the SGH method, and by an exact/approximate Riemann solver in the CCH method. Although these two methods have been developed over decades, this paper will provide a general and systematic scheme that can be applied to both the SGH and CCH methods without any empirical parameters or Riemann solvers.

This scheme contains two key points: one is the relationship between pressure and velocity. Besides the fact that the relationship (\ref{pptautau}) between pressure and specific volume has been constructed, the first key point requires the connection between specific volume and velocity. The other is known as Newton's second law, which is applied to calculate the velocity. Therefore, the new scheme contains the following three steps:
\begin{itemize}
\item Step 1: Construct the relationship between pressure and velocity by using specific volume as an intermediary.
\item Step 2: Apply Newton's second law to calculate the velocity.
\item Step 3: Update the gas dynamic systems in both SGH and CCH methods.
\end{itemize}

\begin{figure}[!h]
\centering
	\includegraphics[scale=0.8]{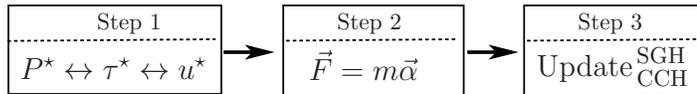}
	\caption{Three steps for SGH and CCH methods by using the same scheme: 1) Build the relationship between pressure and velocity by specific volume; 2) Apply Newton's second law to obtain the velocity; 3) Update the corresponding equations for SGH and CCH methods.\label{stagcent}}	
\end{figure}

To unify both SGH and CCH methods, the strategy for constructing the scheme can be formulated as shown in Figure \ref{stagcent}, where $\vec{F}$ is the force and $\vec{\alpha}$ is the acceleration. 

\subsection{SGH method}
The SGH method may utilize the above three strategy and follow the construction of the scheme while introducing the elements of the new approach. 

\textbf{Step 1.} According to the motion of the nodes in the characteristic time $\delta{t}_{j+\frac{1}{2}}$, the change of specific volume of cell $\Omega_{j+\frac{1}{2}}$ can be expressed as
\begin{eqnarray}
\label{tau-tau0}
\tau^{\star}_{j+\frac{1}{2}}  - \tau_{j+\frac{1}{2}}
 &=& \frac{ V^{\star}_{j+\frac{1}{2}}  - V_{j+\frac{1}{2}} } {m_{j+\frac{1}{2}}}   \nonumber \\
 &=& \frac{(u_{j+1}\textbf{n}^{j+1}_{j+\frac{1}{2}}+u_{j}\textbf{n}^j_{j+\frac{1}{2}})\delta{t}_{j+\frac{1}{2}}}{m_{j+\frac{1}{2}}}  \quad {\rm due \ to \ (\ref{np1})}   \nonumber\\
 &=& \frac{(u_{j+1}-u_{j})\delta{t}_{j+\frac{1}{2}}}{m_{j+\frac{1}{2}}}  ~\quad\quad\quad\quad\quad {\rm owing \ to  \ (\ref{np11})}   \nonumber\\
 &=& \frac{(u_{j+1} - u_{j})(x_{j+1}-x_j)}{m_{j+\frac{1}{2}}c_{j+\frac{1}{2}}}  \quad\quad\quad {\rm thanks \ to\  (\ref{charatime})}   \nonumber\\
 &=&  \frac{u_{j+1}-u_{j}}{\rho_{j+\frac{1}{2}}c_{j+\frac{1}{2}}}.
\end{eqnarray}
Substituting (\ref{tau-tau0}) into (\ref{pptautau}) reads
\begin{eqnarray}
\label{ptau2ndu}
P_{j+\frac{1}{2}}^{\star}=P_{j+\frac{1}{2}} -  \rho_{j+\frac{1}{2}}c_{j+\frac{1}{2}}({u_{j+1}-u_{j}}) + \frac{\gamma+1}{2}\rho_{j+\frac{1}{2}}({u_{j+1}-u_{j}})^2.
\end{eqnarray}

In the case of a cell undergoing an isentropic process, characterized by reversible thermodynamical processes, such as rarefaction waves or isentropic compression, the discretization of the internal energy equation needs to satisfy the entropy conservation. This allows to set $P_{j+\frac{1}{2}}^{\star} = P_{j+\frac{1}{2}}$, which will cancel the entropy production in the case of smooth flows. Finally, the relationship between pressure and velocity writes
\begin{eqnarray}
\label{pfinal}
P_{j+\frac{1}{2}}^{\star}=\left\{
\begin{array}{lll}
P_{j+\frac{1}{2}} -  \rho_{j+\frac{1}{2}}c_{j+\frac{1}{2}}({u_{j+1}-u_{j}}) &{\rm if } \ {u_{j+1}-u_{j}}<0,  \\
~~~~~~ + \frac{\gamma+1}{2}\rho_{j+\frac{1}{2}}({u_{j+1}-u_{j}})^2, \\
P_{j+\frac{1}{2}}   &{\rm if } \ {u_{j+1}-u_{j}}\geq0.
\end{array}
\right.
\end{eqnarray}

The above result (\ref{pfinal}) recovers the Kuropatenko artificial viscosity \cite{kuropatenko_1967}, who derived the formula for computing the pressure jump produced by only one shock wave created by a velocity jump $u_{j+1}-u_{j}$. More sophisticated artificial viscosity models are available (as described, e.g., by Campbell and Shashkov [6]), the simple linear-plus-quadratic model is sufficient to demonstrate the efficacy of the numerical methods for the 1-D gas dynamics problems discussed here.

\textbf{Step 2.} Applying Newton's second law to the control volume $[x_{j-\frac{1}{2}},x_{j+\frac{1}{2}}]$ of node $x_j$ yields
\begin{eqnarray}
\label{slsdis1m}
\alpha_{j}  = \frac{\vec{F}_j}{m_j} =  -\frac{1}{m_j} \left( {P_{j+\frac{1}{2}}^{\star}\textbf{n}_{j+\frac{1}{2},l}^{j} + P^{\star}_{j-\frac{1}{2}}}\textbf{n}_{j-\frac{1}{2},r}^{j+1}\right).
\end{eqnarray}
This allows to get the average velocity of a vertex $x_j$ as
\begin{eqnarray}
\label{slsdis2m}
u^{\star}_j = u^n_j + \frac{1}{2}\alpha_{j}\Delta{t}. 
\end{eqnarray}

\textbf{Step 3.} The discretization of the gas dynamics system can be updated by using (\ref{slsdis1m}) and (\ref{slsdis2m}). Generally, the SGH method employs a predictor-corrector method in time. Knowing all physical quantities at time $t^n$, the resulting numerical scheme can be written as 
\begin{framed}
\center \textbf{Predictor step}
\begin{eqnarray}
\label{slsdis2a}
u^{n+\frac{1}{2}}_{j} &=&2 u^{\star}_{j}  - u^{n}_{j}, \\
\label{slsdis2}
\varepsilon^{n+\frac{1}{2}}_{j+\frac{1}{2}} &=& \varepsilon^{n}_{j+\frac{1}{2}} - \frac{\Delta{t}}{m_{j+\frac{1}{2}}}P_{j+\frac{1}{2}}^{\star}\left(u_{j}^{\star}\textbf{n}_{j+\frac{1}{2}}^{j}
+ u_{j+1}^{\star}\textbf{n}_{j+\frac{1}{2}}^{j+1}\right), \\
\label{slsdis3}
x^{n+\frac{1}{2}}_{j} &=& x^n_j +  u_{j+1}^{\star}\Delta{t}, \\
\label{slsdis4}
\rho^{n+\frac{1}{2}}_{j+\frac{1}{2}} &=&  \frac{\rho^{n}_{j+\frac{1}{2}} V^{n}_{j+\frac{1}{2}} }{V^{n+\frac{1}{2}}_{j+\frac{1}{2}} }, \\
\label{slsdis5}
P^{n+\frac{1}{2}}_{j+\frac{1}{2}} &=& {\rm EOS}\left( \rho^{n+\frac{1}{2}}_{j+\frac{1}{2}} ,\varepsilon^{n+\frac{1}{2}}_{j+\frac{1}{2}} \right).
\end{eqnarray}
\end{framed}

Then applying the equations (\ref{pfinal}) and (\ref{slsdis2m}) again with $P^{n+\frac{1}{2}},\rho^{n+\frac{1}{2}},c^{n+\frac{1}{2}},u^{n+\frac{1}{2}}$ at time $t^{n+\frac{1}{2}}$, the corrector step is written as 
\begin{framed}
\center \textbf{Corrector step}
\begin{eqnarray}
\label{aslsdis0}
u^{n+1}_{j} &=& {2} u^{\ast}_{j} -  u^{n}_{j} , \\
\label{aslsdis2}
\varepsilon^{n+1}_{j+\frac{1}{2}} &=& \varepsilon^{n}_{j+\frac{1}{2}} - \frac{\Delta{t}}{m_{j+\frac{1}{2}}}\frac{P_{j+\frac{1}{2}}^{\star}+P_{j+\frac{1}{2}}^{\ast}}{2}\left(u_{j}^{\ast}\textbf{n}_{j+\frac{1}{2}}^{j}+ u_{j+1}^{\ast}\textbf{n}_{j+\frac{1}{2}}^{j+1}\right), \\
\label{aslsdis3}
x^{n+1}_{j} &=& x^n_j +  u_{j+1}^{\ast}\Delta{t}, \\
\label{aslsdis4}
\rho^{n+1}_{j+\frac{1}{2}} &=&  \frac{\rho^{n}_{j+\frac{1}{2}} V^{n+1}_{j+\frac{1}{2}} }{V^{n}_{j+\frac{1}{2}} }, \\
\label{aslsdis5}
P^{n+1}_{j+\frac{1}{2}} &=& {\rm EOS}\left( \rho^{n+1}_{j+\frac{1}{2}} ,\varepsilon^{n+1}_{j+\frac{1}{2}} \right).
\end{eqnarray}
\end{framed}

Three steps of the new scheme have been presented above for the SGH method. An important observation is that the “artificial viscosity term” no longer contains the empirical parameters, such as $c_1$ and $c_2$ in \cite{caramana_formulations_1998}. This "artificial viscosity term" comes from the Taylor expansion of the state ($\tau,P$) and a link between specific volume and velocity. 
 
\subsection{CCH method} The same strategy mentioned above will be applied to the CCH method.

\textbf{Step 1.}
The relationship between specific volume and velocity will be constructed to obtain the link between pressure and velocity. Suppose a perturbation is created by the collision of two adjacent sub-cells $\Omega_{j-\frac{1}{2},r}$ and $\Omega_{j+\frac{1}{2},l}$ located at vertex $x_j$. In the character time $\delta{t}_{j+\frac{1}{2},l}$, the velocity of the left endpoint of the sub-cell $\Omega_{j+\frac{1}{2},l}$ changes from $u_{j+\frac{1}{2}}$ to $u^{\star}_j$, while the right endpoint remains unchanged as $u_{j+\frac{1}{2}}$. This allows to build the relationship between the specific volume and nodal velocity with the equation (\ref{charatime}) as
\begin{eqnarray}
\label{c1tau-tau0}
\tau^{\star}_{j+\frac{1}{2},l}  - \tau_{j+\frac{1}{2},l} 
&=& \frac{V^{\star}_{j+\frac{1}{2},l}  - V_{j+\frac{1}{2},l} }{m_{j+\frac{1}{2},l}} \nonumber  \\
&=&\frac{(u^{\star}_j-{u_{j+\frac{1}{2}}} )\textbf{n}^j_{{j+\frac{1}{2}},l}\delta{t}_{j+\frac{1}{2},l} }{m_{j+\frac{1}{2},l}} \quad\quad {\rm using\ (\ref{np11}) }\nonumber \\
&=&\frac{(-u^{\star}_j+{u_{j+\frac{1}{2}}} )\delta{t}_{j+\frac{1}{2},l} }{m_{j+\frac{1}{2},l}} \quad\quad\quad\quad {\rm thanks \ to \ (\ref{charatime})}\nonumber \\
&= &\frac{(-u^{\star}_j+{u_{j+\frac{1}{2}}})(x_{j+\frac{1}{2}}-x_j)}{\rho_{j+\frac{1}{2}}(x_{j+\frac{1}{2}}-x_j)c_{j+\frac{1}{2}}}\nonumber \\
&= &\frac{-u^{\star}_j+{u_{j+\frac{1}{2}}} }{\rho_{j+\frac{1}{2}}c_{j+\frac{1}{2}}}
\end{eqnarray}
Similarly, for the sub-cell $\Omega_{j-\frac{1}{2},r}$, one gets
\begin{eqnarray}
\label{c2tau-tau0}
\tau^{\star}_{j-\frac{1}{2},r}  - \tau_{j-\frac{1}{2},r} =  \frac{(u^{\star}_j-{u_{j-\frac{1}{2},r}} )\textbf{n}^j_{{j-\frac{1}{2}},r}\delta{t}_{j-\frac{1}{2},r} }{m_{j-\frac{1}{2},r}} =  \frac{u^{\star}_j-{u_{j-\frac{1}{2}}}}{\rho_{j-\frac{1}{2}}c_{j-\frac{1}{2}}}.
\end{eqnarray}

Substituting (\ref{c1tau-tau0}) and (\ref{c2tau-tau0}) into (\ref{pptautau}) and taking the first-order approximation yields
\begin{eqnarray}
\label{pstarua1}
 P^{\star,1st}_{j+\frac{1}{2},l} =  P_{j+\frac{1}{2}} + \rho_{j+\frac{1}{2}} c_{j+\frac{1}{2}}\left(u^{\star,1st}_j-{u_{j+\frac{1}{2}}}\right) , \\
\label{pstarua2}
P^{\star,1st}_{j-\frac{1}{2},r} =  P_{j+\frac{1}{2}} - \rho_{j+\frac{1}{2}} c_{j+\frac{1}{2}} \left(u^{\star,1st}_j-{u_{j-\frac{1}{2}}}\right).
\end{eqnarray}

Substituting (\ref{c1tau-tau0}) and (\ref{c2tau-tau0}) into (\ref{pptautau}) and taking the first- and second-order approximations reads
\begin{eqnarray}
\label{pstaru1}
 P^{\star,2nd}_{j+\frac{1}{2},l} =  P_{j+\frac{1}{2}} + \rho_{j+\frac{1}{2}} c_{j+\frac{1}{2}}(u^{\star,2nd}_j-{u_{j+\frac{1}{2}}})  + \frac{\gamma+1}{2}\rho_{j+\frac{1}{2}}(u^{\star,2nd}_j-{u_{j+\frac{1}{2}}})^2, \\
\label{pstaru2}
P^{\star,2nd}_{j-\frac{1}{2},r} =  P_{j+\frac{1}{2}} - \rho_{j+\frac{1}{2}} c_{j+\frac{1}{2}}(u^{\star,2nd}_j-{u_{j-\frac{1}{2}}}) + \frac{\gamma+1}{2}\rho_{j-\frac{1}{2}}(u^{\star,2nd}_j-{u_{j-\frac{1}{2}}})^2.
\end{eqnarray}

\textbf{Step 2.}
Since the velocity of the vertex $x_j$ is a constant within each time step, this condition implies that the nodal acceleration is $\vec{\alpha}_j = 0$. According to this new finding, Newton's second law is applied to the control volume $[x_{j-\frac{1}{2}},x_{j+\frac{1}{2}}]$ of node $x_j$, one gets
\begin{eqnarray}
\label{pp0}
0 = m_j\vec{\alpha}_j = \vec{F}_j = -\left(P^{\star}_{j-\frac{1}{2},r} \textbf{n}^{j-\frac{1}{2}}_{j-\frac{1}{2},r} +  P^{\star}_{j+\frac{1}{2},l} \textbf{n}^{j+\frac{1}{2}}_{j+\frac{1}{2},l} \right)
\end{eqnarray}

Taking (\ref{pstarua1}) and (\ref{pstarua2}) into (\ref{pp0}) reads
\begin{eqnarray}
\label{pufinal1}
&&u^{\star,1st}_{j}=\frac{ \rho_{j+\frac{1}{2}}c_{j+\frac{1}{2}} +  \rho_{j-\frac{1}{2}}c_{j-\frac{1}{2}} + P_{j-\frac{1}{2}}-P_{j+\frac{1}{2}}  }{\rho_{j+\frac{1}{2}}c_{j+\frac{1}{2}} +  \rho_{j-\frac{1}{2}}c_{j-\frac{1}{2}}}, \\
&&P^{\star,1st}_{j-\frac{1}{2},r}=P^{\star,1st}_{j+\frac{1}{2},l}=\frac{ \rho_{j+\frac{1}{2}}c_{j+\frac{1}{2}}P_{j-\frac{1}{2}} + \rho_{j-\frac{1}{2}}c_{j-\frac{1}{2}}P_{j+\frac{1}{2}} }{\rho_{j+\frac{1}{2}}c_{j+\frac{1}{2}} +  \rho_{j-\frac{1}{2}}c_{j-\frac{1}{2}}} \nonumber \\
&&\quad\quad\quad\quad\quad\quad\quad\quad-\frac{  \rho_{j+\frac{1}{2}}c_{j+\frac{1}{2}}\rho_{j-\frac{1}{2}}c_{j-\frac{1}{2}}\left( u_{j+\frac{1}{2}} - u_{j-\frac{1}{2}} \right)}{\rho_{j+\frac{1}{2}}c_{j+\frac{1}{2}} +  \rho_{j-\frac{1}{2}}c_{j-\frac{1}{2}}}.\label{pufinal2}
\end{eqnarray}
It is noticed that the above equations are the acoustic solver \cite{godunov_1959}. 

Similarly, substituting (\ref{pstaru1}) and (\ref{pstaru2}) into (\ref{pp0}) writes
\begin{eqnarray}
&A\left(u^{\star,2nd}_j\right)^2 + Bu^{\star,2nd}_j + C = 0 \\
&\left\{\begin{array}{llll}
A = \frac{\gamma+1}{2}\left(\rho_{j-\frac{1}{2}} - \rho_{j+\frac{1}{2}}\right), \\
B = -\left[(\gamma+1)(\rho_{j-\frac{1}{2}}u_{j-\frac{1}{2}}-\rho_{j+\frac{1}{2}}u_{j+\frac{1}{2}})
+ \rho_{j-\frac{1}{2}}c_{j-\frac{1}{2}} + \rho_{j+\frac{1}{2}}c_{j+\frac{1}{2}}\right], \\
C = \frac{\gamma+1}{2}\left(\rho_{j-\frac{1}{2}}u^2_{j-\frac{1}{2}}-\rho_{j+\frac{1}{2}}u^2_{j+\frac{1}{2}}\right) + P_{j-\frac{1}{2}} - P_{j+\frac{1}{2}} \\
\quad \quad + \rho_{j-\frac{1}{2}}u_{j-\frac{1}{2}}c_{j-\frac{1}{2}}+\rho_{j+\frac{1}{2}}u_{j+\frac{1}{2}}c_{j-\frac{1}{2}} .
\end{array}\right.\nonumber
\end{eqnarray}
The above equation has two solutions
\begin{eqnarray}
\label{ustar2s}
u^{\star,2nd}_{j_{\pm}}  = \frac{-B\pm\sqrt{B^2-4AC}}{2A}, \quad {\rm if } ~ B^2-4AC>0 ~ {\rm and} ~ A\neq0.
\end{eqnarray}
A suitable solution will be chosen as follows
 \begin{eqnarray}
 \label{ustar2nd}
|u^{\star,2nd}_j-u^{1st}_j|=\min\left(|u^{\star,2nd}_{j_{+}}-u^{\star,1st}_j|,|u^{\star,2nd}_{j_{-}}-u^{\star,1st}_j|\right), \\
 \left\{
\begin{array}{lll}
\rho_{j+\frac{1}{2}} c_{j+\frac{1}{2}} (u^{\star,2nd}_{j}-{u_{j+\frac{1}{2}}})^2  \geq \frac{\gamma+1}{2}\rho_{j+\frac{1}{2}}|u^{\star,2nd}_{j}-{u_{j+\frac{1}{2}}}|^3,\\
\rho_{j+\frac{1}{2}} c_{j+\frac{1}{2}} (u^{\star,2nd}_{j}-{u_{j-\frac{1}{2}}})^2  \geq \frac{\gamma+1}{2}\rho_{j-\frac{1}{2}}|u^{\star,2nd}_{j}-{u_{j-\frac{1}{2}}}|^3. \nonumber
\end{array}
\right.
\end{eqnarray}

Once the nodal velocity $u^{\star,2nd}$ is known, the pressures can be calculated using (\ref{pstaru1}) and (\ref{pstaru2}). If $u^{\star,2nd}$ does not satisfy the equations (\ref{ustar2nd}), the velocity and pressures will be replaced by the equations (\ref{pufinal1}) and (\ref{pufinal2}). For the sake of brevity, the nodal velocity and pressures are defined as $u^{\star}_j$ and $P^{\star}_{j+\frac{1}{2},l}(P^{\star}_{j-\frac{1}{2},r})$, respectively. 

\textbf{Step 3.}  The discretization of the gas dynamics system can be updated by using $u^{\star}_j$ and $P^{\star}_{j+\frac{1}{2},l}(P^{\star}_{j+\frac{1}{2},r})$ as follows 
\begin{framed}
\center \textbf{Forward Euler step}
\begin{eqnarray}
\label{clsdis1}
u^{n+1}_{j+\frac{1}{2}} &=& u^{n}_{j+\frac{1}{2}} - \frac{\Delta{t}}{m_{j+\frac{1}{2}}}\left({P_{j+\frac{1}{2},l}^{\star}\textbf{n}_{j +\frac{1}{2},l}^{j}
+ P^{\star}_{j+\frac{1}{2},r}}\textbf{n}_{j+\frac{1}{2},r}^{j+1} \right), \\
\label{clsdis2}
E^{n+1}_{j+\frac{1}{2}} &=& E^{n}_{j+\frac{1}{2}} - \frac{\Delta{t}}{m_{j+\frac{1}{2}}}\left({P_{j+\frac{1}{2},l}^{\star}u_{j}^{\star} \textbf{n}_{j+\frac{1}{2},l}^{j}
+P^{\star}_{j+\frac{1}{2},r}u_{j+1}^{\star}}\textbf{n}_{j+\frac{1}{2},r}^{j+1}\right), \\
\label{clsdis3}
x^{n+1}_{j} &=& x^n_j +  u_{j+1}^{\star}\Delta{t}, \\
\label{clsdis4}
\rho^{n+1}_{j+\frac{1}{2}} &=&  \frac{\rho^{n}_{j+\frac{1}{2}} V^{n+1}_{j+\frac{1}{2}} }{V^{n}_{j+\frac{1}{2}} }, \\
\label{clsdis5}
P^{n+1}_{j+\frac{1}{2}} &=& {\rm EOS}\left( \rho^{n+1}_{j+\frac{1}{2}} ,\varepsilon^{n+1}_{j+\frac{1}{2}} \right).
\end{eqnarray}
\end{framed}

The new finding that the nodal acceleration is zero in each time step is employed for the CCH method. Instead of an exact/approximate Riemann solver, Newton's second law can be applied to calculate the velocity and pressure. Furthermore, this form can be easily extended to multi-dimensions, which constitutes a certain departure from the standard dimension-by-dimension approach.

In summary, it is important to emphasize here that the same scheme constructed above clearly unifies both SGH and CCH methods. The scheme neither contains empirical parameters nor solves the Riemann problem, which is the main contribution of this work.

\subsection {Main properties}
Following the same strategy developed when constructing both the SGH and CCH methods, the two following tasks will be performed 
\begin{itemize}
\item Conservatively.
\item Compatibility with the entropy condition.
\end{itemize}

Following this path, the essential properties of both the SGH and CCH methods will be verified, such as the conservation of total mass, momentum, and energy. For the sake of brevity, the predictor step will be checked in the SGH method since the corrector step has the same properties.

\subsubsection{Conservation of total mass} This is obtained directly from the equations (\ref{slsdis4})  (\ref{aslsdis4}) and (\ref{clsdis4}).

\subsubsection{Conservation of total momentum} For the SGH method, considering the equations (\ref{slsdis1m}), (\ref{slsdis2m}) and (\ref{slsdis2a}) and taking the sum of all dual cells reads
\begin{eqnarray}
\label{sghmom}
\sum_j^Nm_ju^{n+1}_j = \sum_j^Nm_ju^{n}_j - \sum_j^N(P^{\star}_{j-\frac{1}{2}}-P^{\star}_{j+\frac{1}{2}}) = \sum_j^Nm_ju^{n}_j - {\rm BND.}
\end{eqnarray}
Similarly for the CCH method, summing up all the primary cells with (\ref{clsdis1}) yields
\begin{eqnarray}\label{cchmom}
\sum_j^Nm_{j+\frac{1}{2}}u^{n+1}_{j+\frac{1}{2}} &= &\sum_j^Nm_ju^{n}_{j+\frac{1}{2}} - \sum_j^N\left(P^{\star}_{j+\frac{1}{2},l}-P^{\star}_{j+\frac{1}{2},r}\right)  \\
&=& \sum_j^Nm_ju^{n}_{j+\frac{1}{2}} - {\rm BND.} \nonumber
\end{eqnarray}

\subsubsection{Conservation of total energy} According to the definition of energy for the SGH method, the total energy of all cells can be rewritten as
\begin{eqnarray}
&&\sum_j^N\left(m_{j+\frac{1}{2}}e^{n+1}_{j+\frac{1}{2}} + \frac{1}{2}m_j(u^{n+1}_{j})^2\right) \quad\quad\quad\quad {\rm thanks \ to (\ref{slsdis2})\ and\ (\ref{sghmom})}\nonumber  \\
&&=\sum_j^N\left(m_{j+\frac{1}{2}}e^{n}_{j+\frac{1}{2}} + \frac{1}{2}m_j(u^{n}_{j})^2\right) - \sum_j^N\left[P^{\star}_{j+\frac{1}{2}}\left(u^{\star}_j-u^{\frac{1}{2}}_{j+1}\right)+ u^{\star}_j\left(P^{\star}_{j-\frac{1}{2}}-P^{\star}_{j+\frac{1}{2}}\right)\right] \nonumber\\
&&=\sum_j^N\left(m_{j+\frac{1}{2}}e^{n}_{j+\frac{1}{2}} + \frac{1}{2}m_j(u^{n}_{j})^2\right) - {\rm BND.} \nonumber
\end{eqnarray}
Similarly, for the CCH method, considering the equation (\ref{clsdis2}) yields
\begin{eqnarray}
&&\sum^N_jm_{j+\frac{1}{2}}E^{n+1}_{j+\frac{1}{2}}  = \sum^N_jm_{j+\frac{1}{2}}E^n_{j+\frac{1}{2}} - \sum^N_j\left(P^{\star}_{j+\frac{1}{2},l}U^{\star}_j-
 P^{\star}_{j+\frac{1}{2},r}U^{\star}_{j+1}\right) \nonumber \\
&&=\sum^N_jm_{j+\frac{1}{2}}E^n_{j+\frac{1}{2}} - {\rm BND}. \nonumber
\end{eqnarray}

Without taking into account boundary conditions (BND), both SGH and CCH methods using the same scheme are conservative in total mass, total momentum, and total energy. 

\subsubsection{Entropy inequality} Now, the compatibility of both SGH and CCH methods with the entropy condition will be studied, which may be related to the stability of the two methods.

Using the Gibbs formula for the SCH method, the time variation of the specific entropy can be expressed as 
\begin{eqnarray}
\label{sghentropy}
&&T_{j+\frac{1}{2}}\frac{d}{dt}S_{j+\frac{1}{2}} = \frac{d}{dt}\varepsilon_{j+\frac{1}{2}}  + P_{j+\frac{1}{2}}\frac{d}{dt}\tau_{j+\frac{1}{2}} \\
&&=\frac{d}{dt}\varepsilon_{j+\frac{1}{2}} + P_{j+\frac{1}{2}}\frac{d}{dt}\tau_{j+\frac{1}{2}}
= \left(P_{j+\frac{1}{2}} - P^{\star}_{j+\frac{1}{2}}\right)\left(u_{j+1}-u_j\right)  \nonumber \\
&&=\left\{
\begin{array}{lll}
\rho_{j+\frac{1}{2}}c_{j+\frac{1}{2}}({u_{j+1}-u_{j}})^2 - \frac{\gamma+1}{2}\rho_{j+\frac{1}{2}}({u_{j+1}-u_{j}})^3  &{\rm if }\ {u_{j+1}-u_{j}}<0 \\
0   &{\rm if }\ {u_{j+1}-u_{j}}\geq 0
\end{array} \right.  \nonumber \\
&& \geq  0. \nonumber
\end{eqnarray}
Similarly for the CCH method, the time variation of the specific entropy takes the following form 
\begin{eqnarray}
&&T_{j+\frac{1}{2}}\frac{d}{dt}S_{j+\frac{1}{2}} = \frac{d}{dt}\varepsilon_{j+\frac{1}{2}}  + P_{j+\frac{1}{2}}\frac{d}{dt}\tau_{j+\frac{1}{2}}  \label{cchentropy}\\
&&=  \frac{d}{dt}E_{j+\frac{1}{2}}  - u_{j+\frac{1}{2}}\frac{d}{dt}u_{j+\frac{1}{2}} + P_{j+\frac{1}{2}}\frac{d}{dt}\tau_{j+\frac{1}{2}}  \nonumber\\
&&= \left(P_{j+\frac{1}{2},r} - P_{j+\frac{1}{2}}\right)\left(u^{\star}_{j+1}-u_{j+\frac{1}{2}}\right) +
\left(P_{j+\frac{1}{2},l} - P_{j+\frac{1}{2}}\right)\left(u_{j+\frac{1}{2}}
 - u^{\star}_{j}\right) \nonumber
 \end{eqnarray}
If the $u^{\star}_j$ satisfies the equations (\ref{ustar2nd}), the equation (\ref{cchentropy}) becomes 
 \begin{eqnarray}
 \label{cchentropy1}
T_{j+\frac{1}{2}}\frac{d}{dt}S_{j+\frac{1}{2}}&&= \rho_{j+\frac{1}{2}}c_{j+\frac{1}{2}}\left[\left(u^{\star}_{j+1}-u_{j+\frac{1}{2}}\right)^2+\left(u_{j+\frac{1}{2}} - u^{\star}_{j}\right)^2 \right]  \\
 &&\quad + \frac{\gamma+1}{2}\rho_{j+\frac{1}{2}}\left[\left(u^{\star}_{j+1}-u_{j+\frac{1}{2}}\right)^3+ \left(u_{j+\frac{1}{2}} - u^{\star}_{j}\right)^3 \right] \geq 0.\nonumber
\end{eqnarray}
Otherwise, the equation (\ref{cchentropy}) can be rewritten as
 \begin{eqnarray}
 \label{cchentropy2}
T_{j+\frac{1}{2}}\frac{d}{dt}S_{j+\frac{1}{2}}&&= \rho_{j+\frac{1}{2}}c_{j+\frac{1}{2}}\left[\left(u^{\star}_{j+1}-u_{j+\frac{1}{2}}\right)^2+\left(u_{j+\frac{1}{2}} - u^{\star}_{j}\right)^2 \right] \geq 0.
\end{eqnarray}

The equations (\ref{cchentropy}), (\ref{cchentropy1}) and (\ref{cchentropy2}) show that the scheme satisfies the entropy condition for both SGH and CCH methods. In particular, the SGH method satisfies the entropy conservation in an isentropic process, while the entropy increases in the CCH method.

\section{Numerical examples}
In this section, a suite of challenging test examples are calculated to demonstrate the accuracy and robustness of both SGH and CCH methods by using the new scheme. Although the algorithm supports several types of equations of state, only the ideal gas equation (EOS) is used here
\begin{eqnarray*}
p = (\gamma -1) \rho \varepsilon,\ \
\gamma = \left\{
\begin{array}{l}
5/3,\quad\rm{monoatomic\ gas, } \\
7/5,\quad\rm{diatomic\  gas, }
\end{array}
\right.
\end{eqnarray*}
where $\gamma$ is the adiabatic index. The following examples are the Sod shock tube \cite{sod_survey_1978}, Lax shock tube \cite{lax_weak_1954}, double rarefaction waves \cite{einfeldt_godunov-type_1991}, Sedov blast wave \cite{sedov_similarity_1961},  shock density wave interaction \cite{shuosher-2009},  and LeBlanc shock tube \cite{zhang_positivity-preserving_2010}. As all these examples contain discontinuities or sharp gradients of some kind, this is a natural choice, and a reference solution is also given for use in accurate comparisons. 

Since the presented CCH method is a first-order scheme, the SGH method only employs the 'Predictor Step' in time, such that the numerical errors are only attributed to the spatial discretization. Nevertheless, both methods are enough to get good results. This is probably because the Lagrangian scheme resolves transport almost exactly, even with a first-order method. For some second-order extensions the interested reader refer to \cite{maire_high-order_2009} \cite{cheng_second_2014} \cite{qian_local_2021}.

\subsection {Sod shock tube} 
\begin{figure}[!h]
\centering
	\includegraphics[width=2 in,keepaspectratio]{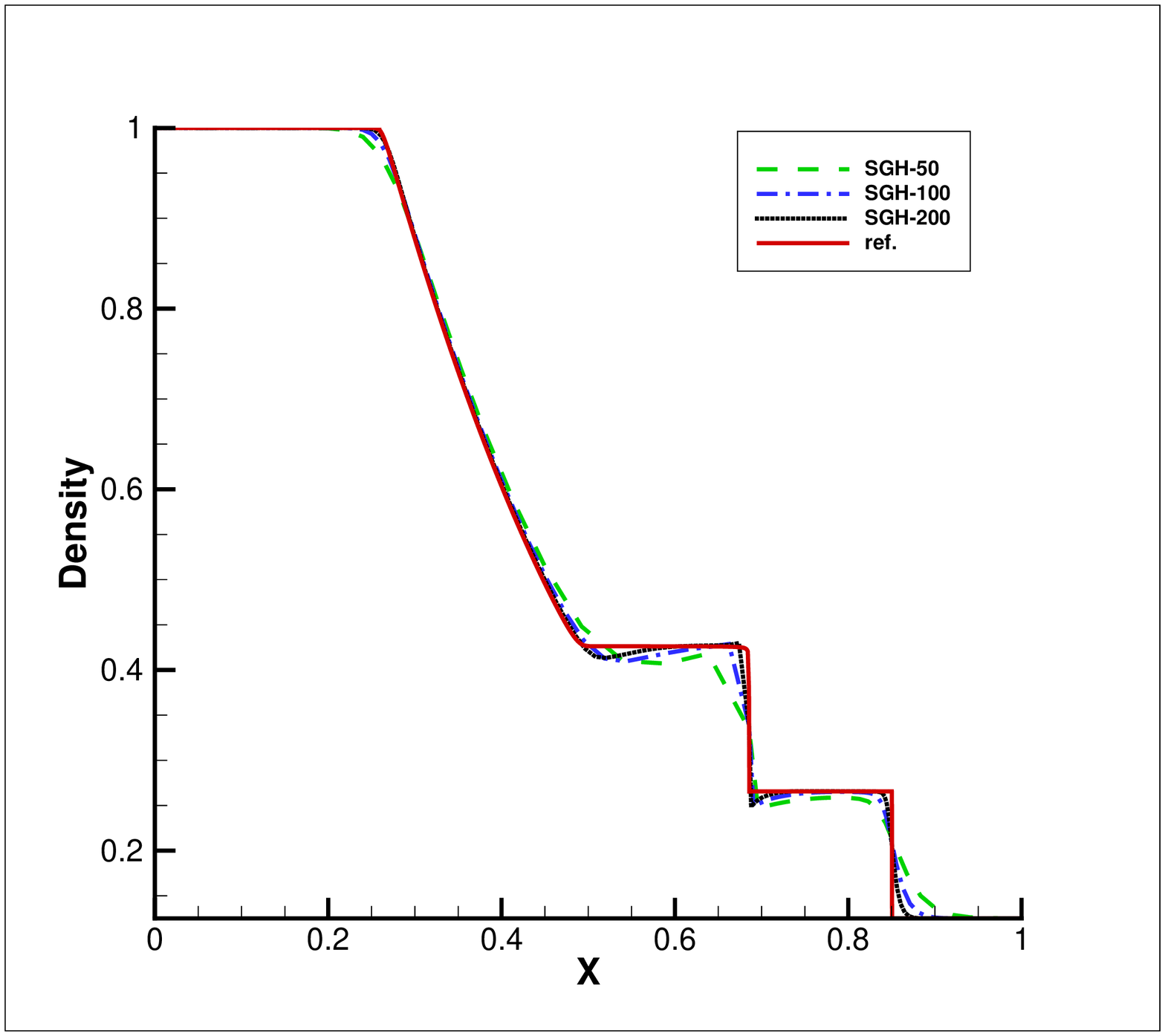}
	\includegraphics[width=2 in,keepaspectratio]{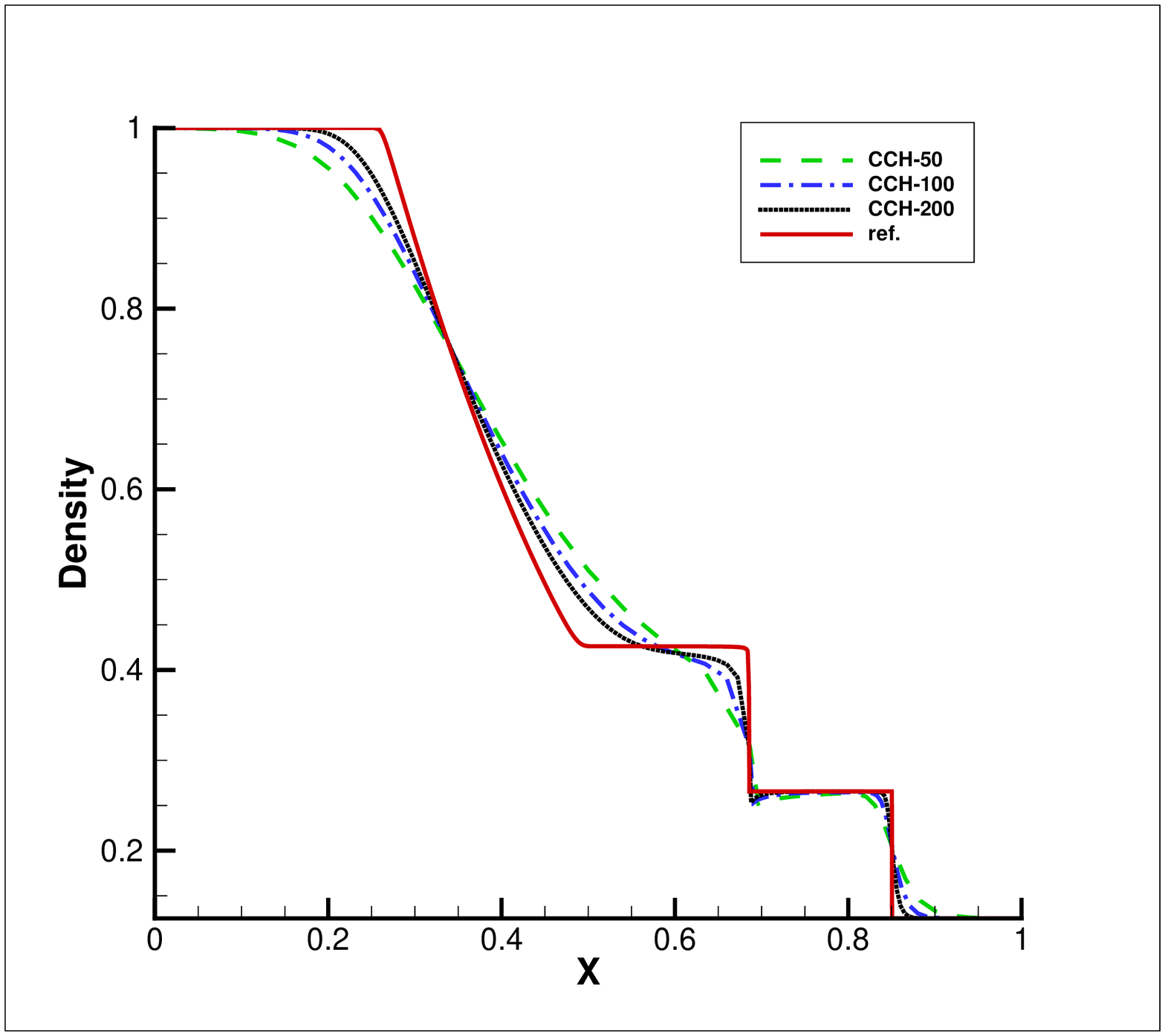}
	\caption{Sod problem: mesh convergence for density profiles with SGH (left) method and CCH (right) method, time $T=0.2$, mesh-size $N=50$ (green\ dashed\ line), $N=100$ (blue\ dashed-dotted\ line), $N=N=200$ (black\ dotted\ line). \label{sod-density-mj}}	
\end{figure}

\begin{figure}[!h]
\centering
	\includegraphics[width=2 in,keepaspectratio]{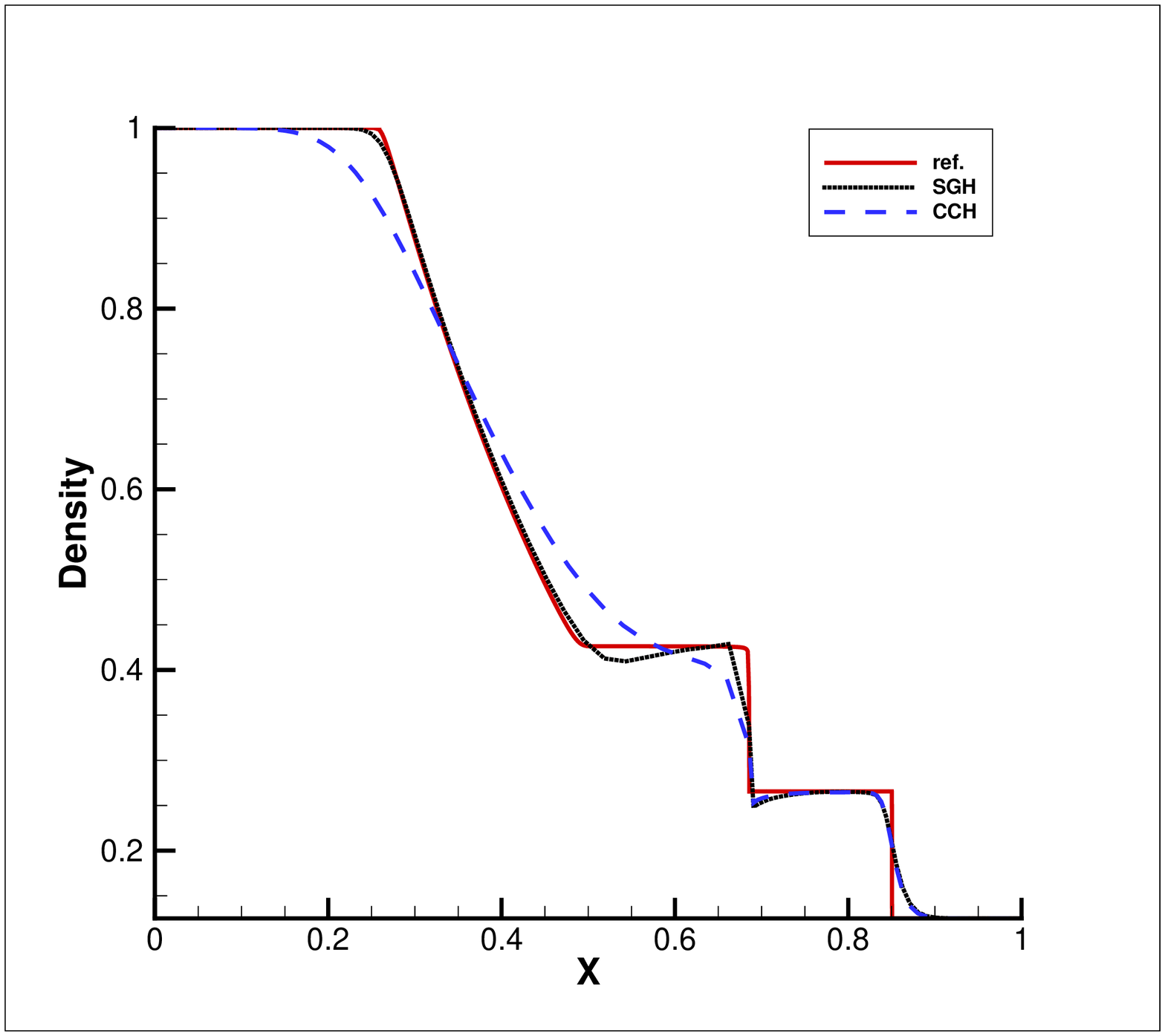}
	\includegraphics[width=2 in,keepaspectratio]{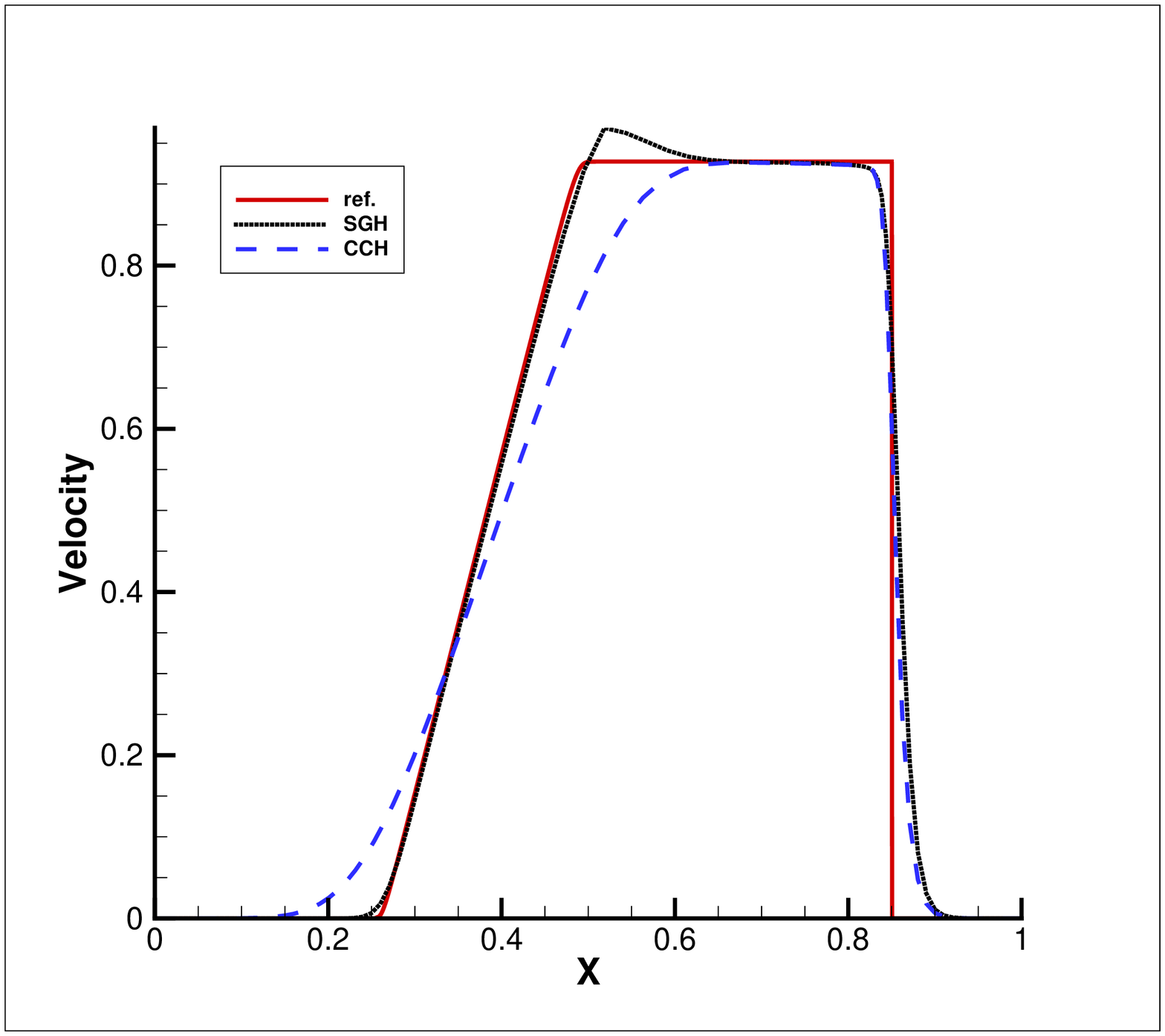}
	\caption{Sod problem: comparison between SGH (black dotted line) method and CCH (blue dashed line) method for density (left) and velocity (right) plots, time $T=0.2$, mesh-size $N=100$. \label{sod-den-vel-vs}}	
\end{figure}

\begin{figure}[!h]
\centering
	\includegraphics[width=2 in,keepaspectratio]{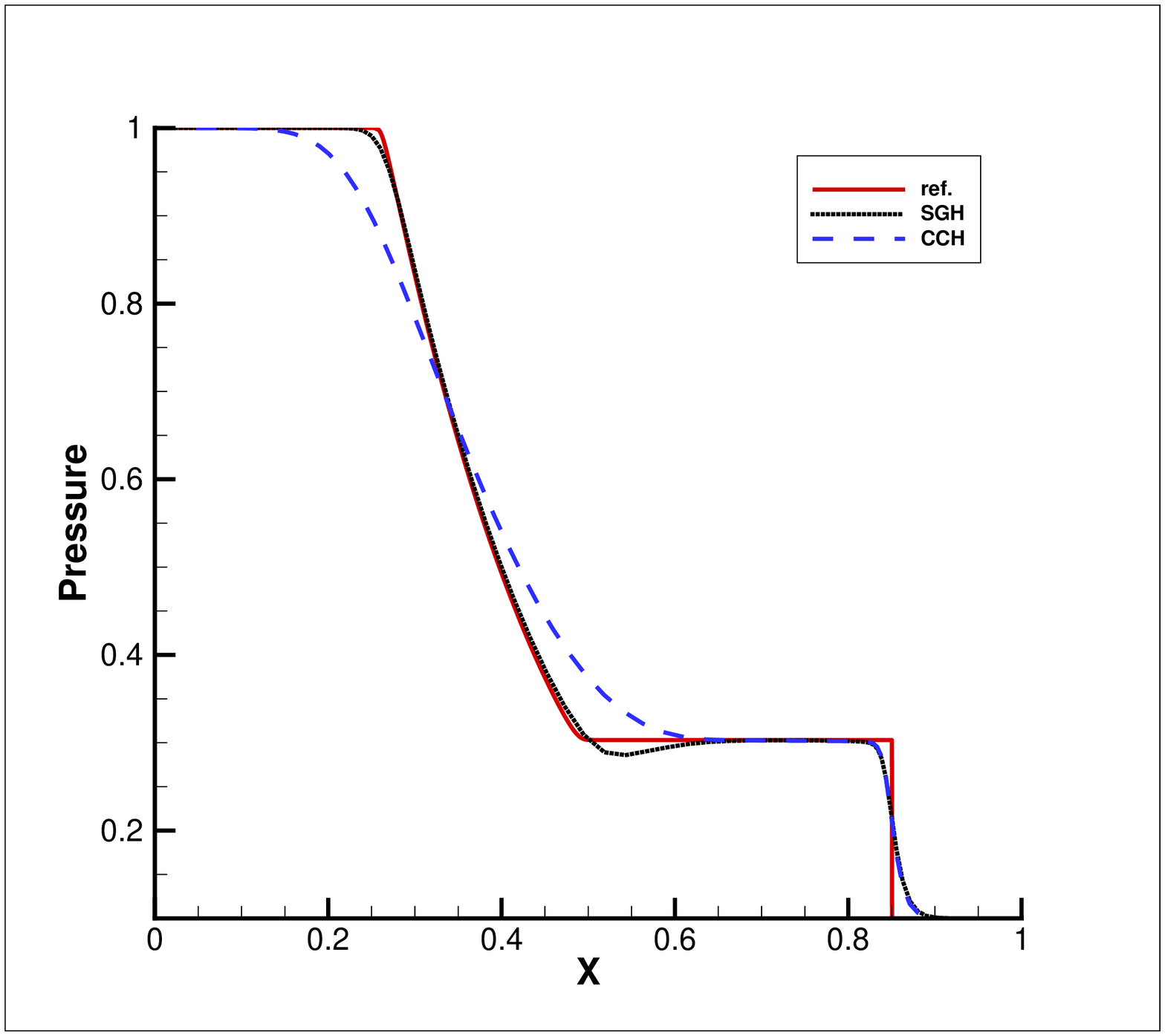}
	\includegraphics[width=2 in,keepaspectratio]{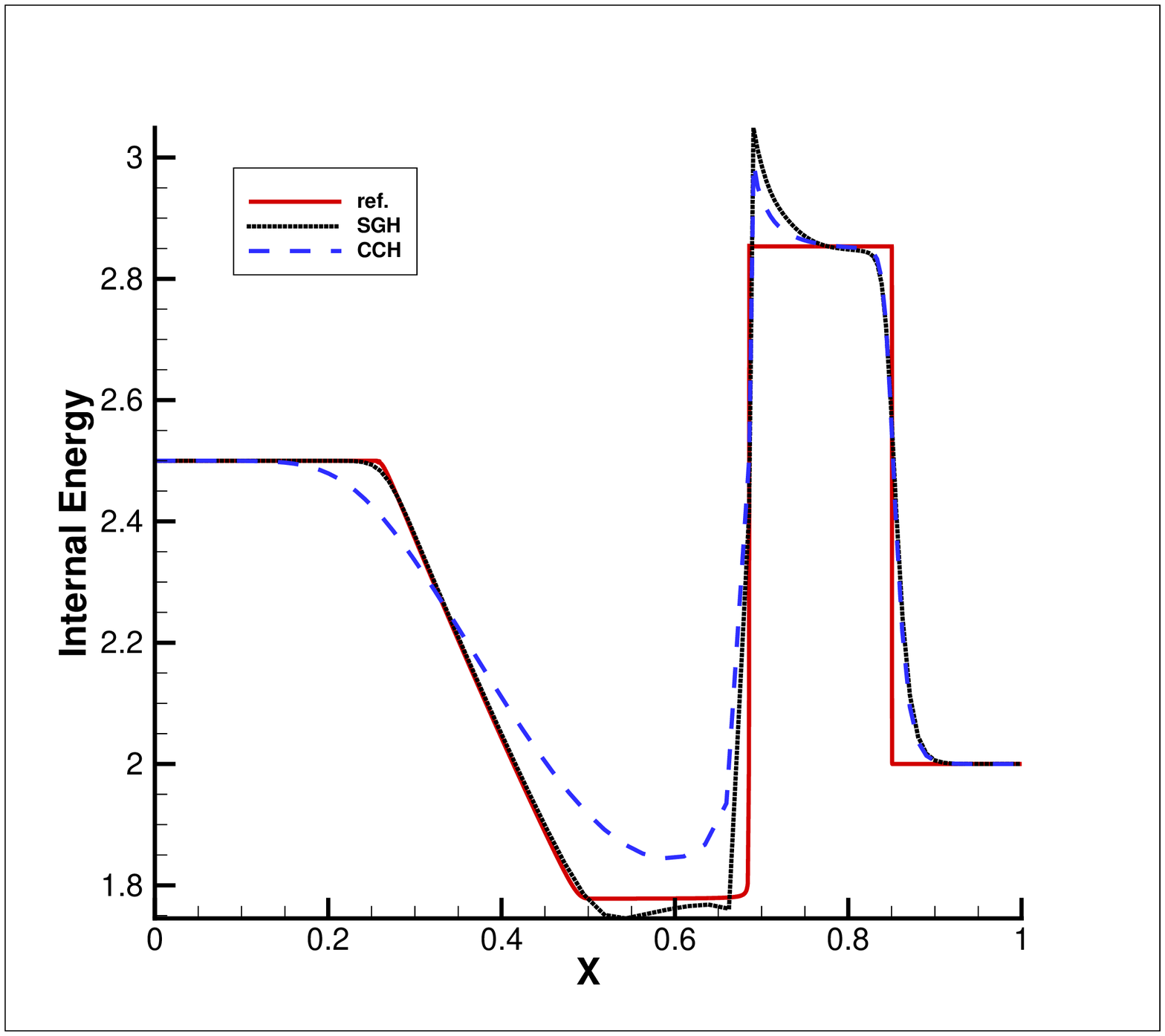}
	\caption{Sod problem: comparison between SGH (black dotted line) method and CCH (blue dashed line) method for pressure (left) and internal energy (right) plots, time $T=0.2$, mesh-size $N=100$. \label{sod-pre-ein-vs}}		
\end{figure}

This is a very mild test with the initial state as 
$$
\begin{array}{ll}
\left(\rho,u,P,\gamma\right)^T=\left\{
\begin{array}{ll}
\left(1,\ 0,\ 1,\ 7/5\right)^T   \quad  &\rm{if} \quad x<0.5,\\
\left(0.125,\ 0,\ 0.1,\ 7/5\right)^T \quad &\rm{if} \quad x\geq0.5.
\end{array}
\right.
\end{array}
$$
(left)
For both SGH and CCH methods, several spatial steps $N=50,100,200$ are presented, which correspond to $N=50,100,200$ with $N$ the number of cells on the unit length in Figure \ref{sod-density-mj}. As expected, when the mesh becomes fine, the numerical results approximate the exact solution. In Figure \ref{sod-den-vel-vs}, compared with the results of the CCH method, the SGH method shows very good accuracy in the rarefaction wave region, even if the two methods are completely based on the same scheme. They also generate slight overshoots near the contact discontinuity in the internal energy plot in Figure \ref{sod-pre-ein-vs}.


\subsection {Lax shock tube}
The one-dimensional Euler equations are simulated with the Riemann initial condition for the Lax problem
$$
\begin{array}{ll}
\left(\rho,u,P,\gamma\right)^T=\left\{
\begin{array}{ll}
\left(0.445,\ 0.698,\ 3.528,\ 7/5\right)^T  \quad &\rm{if} \quad x\leq0.5,\\
\left(0.5,\ 0,\ 0.571,\ 7/5\right)^T \quad &\rm{if} \quad x\geq0.5.
\end{array}
\right.
\end{array}
$$
The solutions are integrated to $T = 0.16$ in the computational domain $[0,1]$. The inlet and outlet conditions are imposed on the left and right boundaries, respectively. The exact solution again comprises a left rarefaction wave, a contact, and a right shock. However, unlike Sod’s test, the contact discontinuity has a rather large jump. 

\begin{figure}[!h]
\centering
	\includegraphics[width=2 in,keepaspectratio]{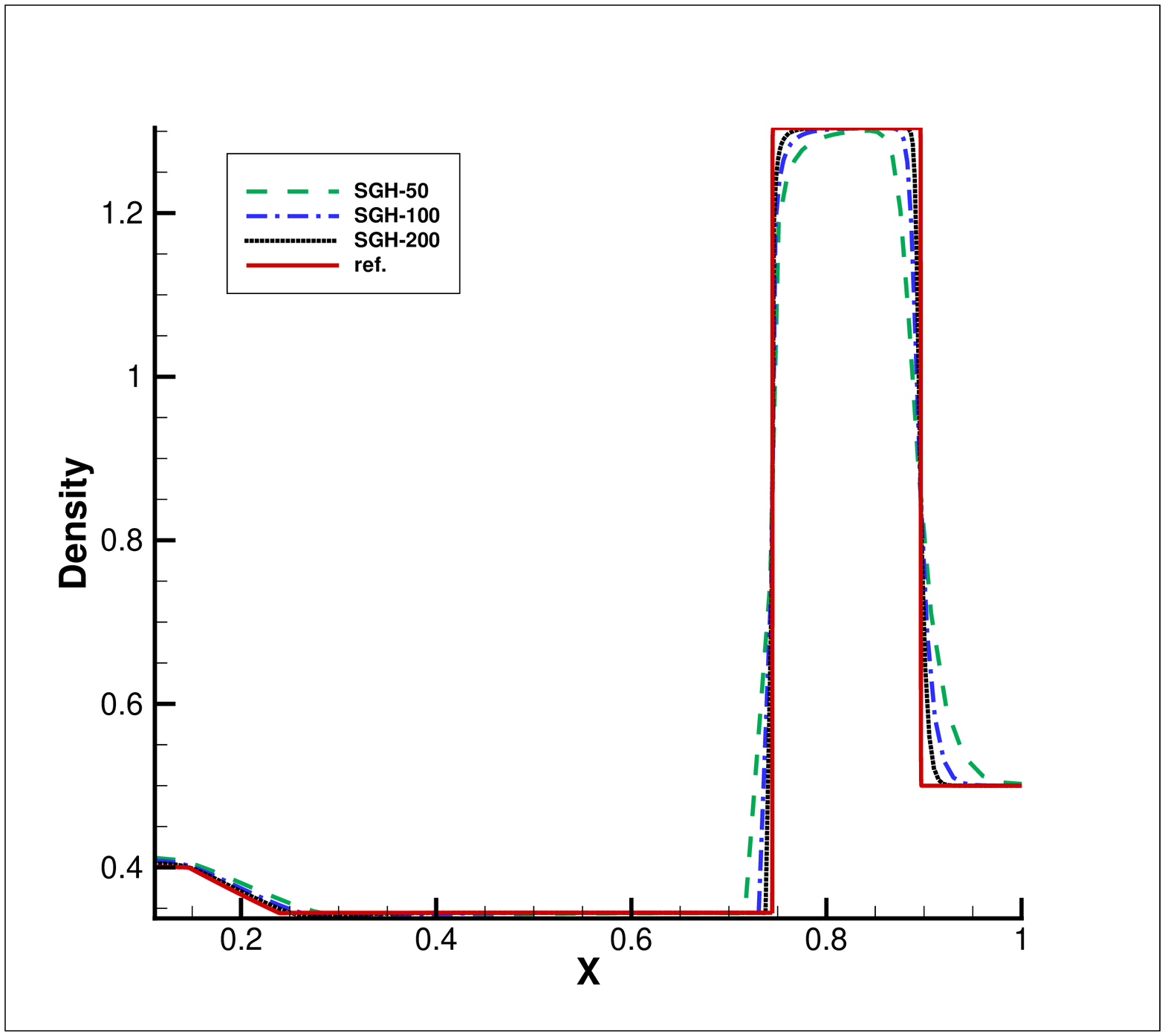}
	\includegraphics[width=2 in,keepaspectratio]{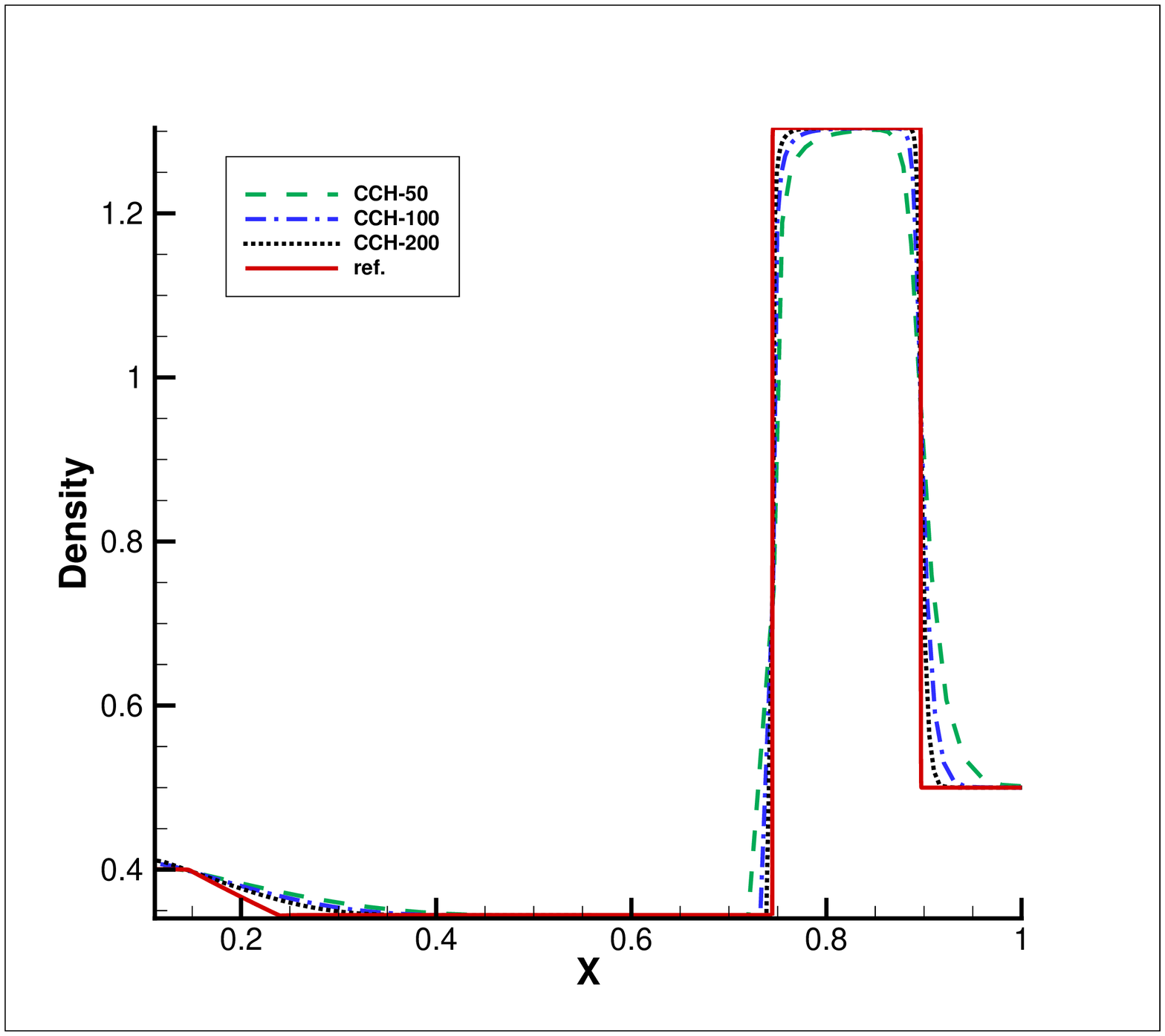}
	\caption{Lax problem: mesh convergence for density profiles with SGH (left) method and CCH (right) method, time $T=0.16$, mesh-size $N=50$ (green\ dashed\ line), $N=100$ (blue\ dashed-dotted\ line), $N=200$ (black\ dotted\ line). \label{lax-density-mj}}		
\end{figure}

\begin{figure}[!h]
\centering
	\includegraphics[width=2 in,keepaspectratio]{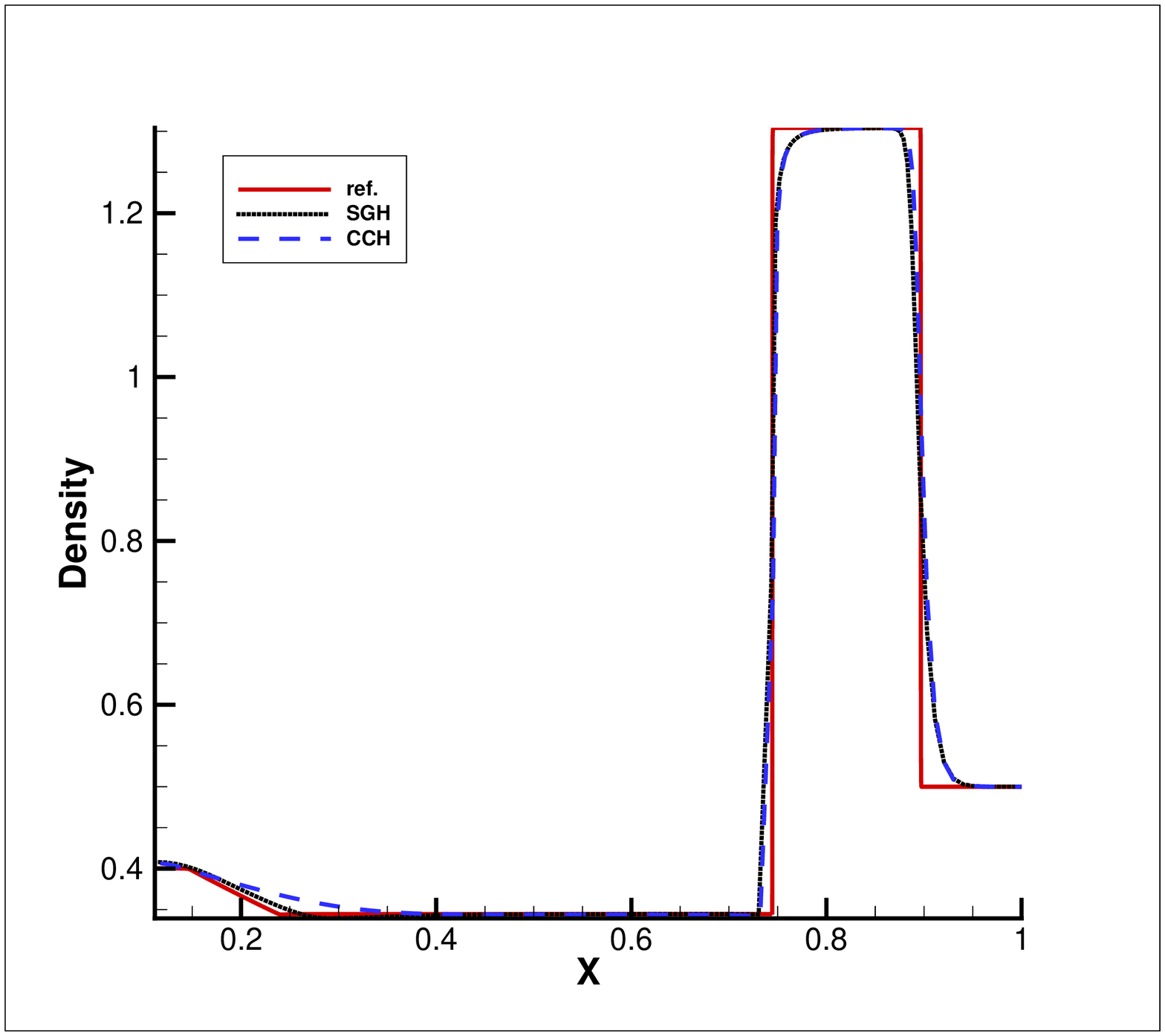}
	\includegraphics[width=2 in,keepaspectratio]{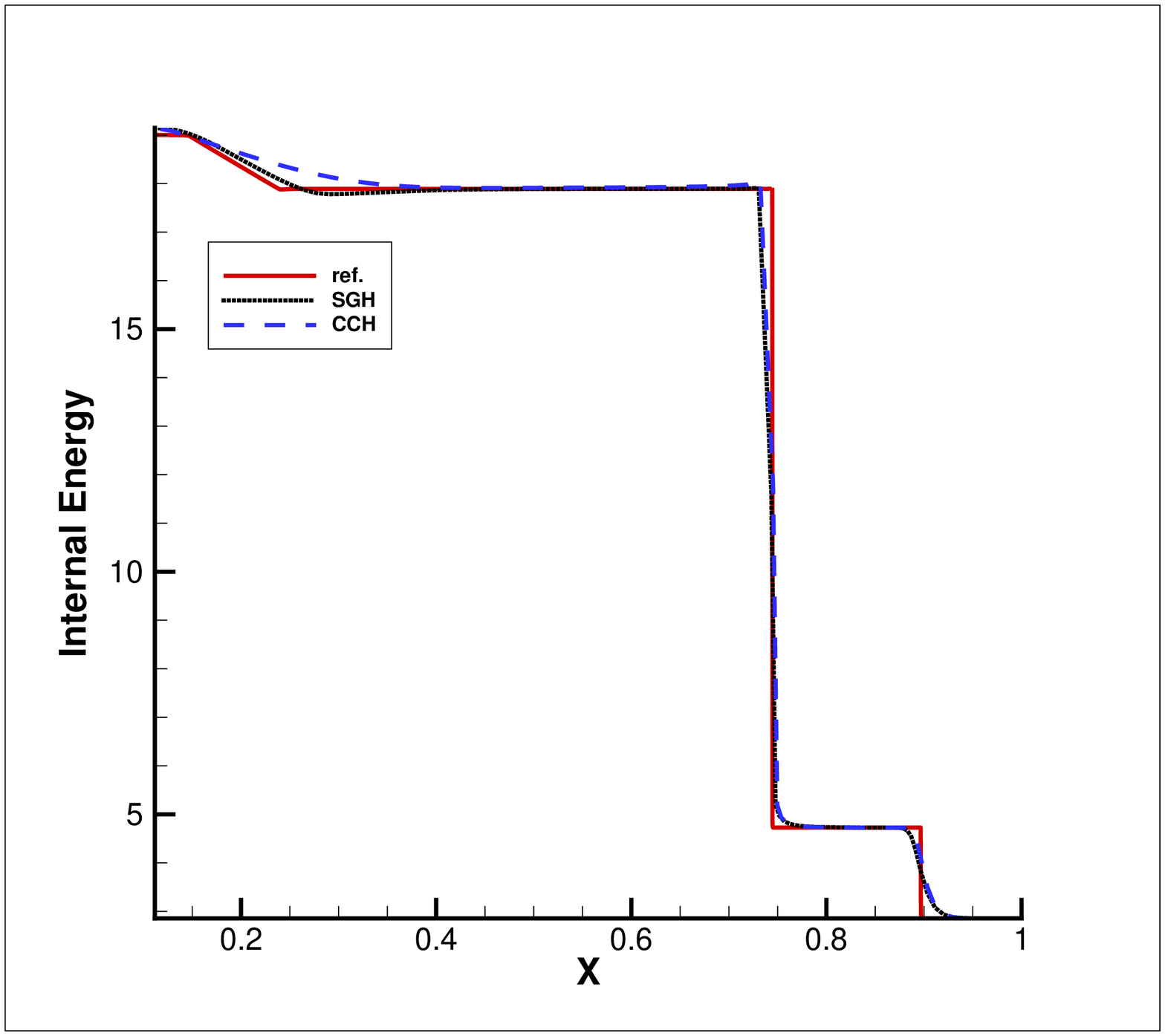}
	\caption{Lax problem: comparison between SGH (black dotted line) method and CCH (blue dashed line) method for density (left) and internal energy (right) plots, time $T=0.16$, mesh-size $N=100$.\label{Lax-pre-ein-vs}}
\end{figure}

Solutions with different grid sizes $N=50,100,200$ are illustrated in Figure \ref{lax-density-mj}. The convergence behavior of both SGH and CCH methods is very satisfactory. In comparison with the numerical solution of the SGH method, the CCH method gives an equivalent solution for the capture of the discontinuities in Figure \ref{Lax-pre-ein-vs}. Similar to the sod problem results, the SGH method provides a much better resolution of left-side rarefaction waves.
\subsection{Double rarefaction waves} 
To check the robustness of the SGH/CCH methods in the presence of highly varying smooth solutions with density and pressure near the vacuum, a classical Riemann problem has two strong rarefaction waves initiated as
$$
\begin{array}{ll}
\left(\rho,u,P,\gamma\right)^T=\left\{
\begin{array}{ll}
\left(1,\ -2,\ 0.4,\ 7/5\right)^T  \quad &\rm{if} \quad x<0.5,\\
\left(1,\ 2,\ 0.4,\ 7/5\right)^T \quad  &\rm{if} \quad x\geq0.5.
\end{array}
\right.
\end{array}
$$
The velocities of the left and right boundaries are prescribed to $u=-2$ and $u = 2$, respectively. The initial conditions are selected to produce vacuum at $x=0.5$. The solution consists of two strong rarefactions with a trivial stationary contact, and the resulting middle state is close to vacuum for the density and pressure. The simulation is implemented in $x\in[0,1]$ at $T = 0.15$.

\begin{figure}[!h]
\centering	
	\includegraphics[width=2 in,keepaspectratio]{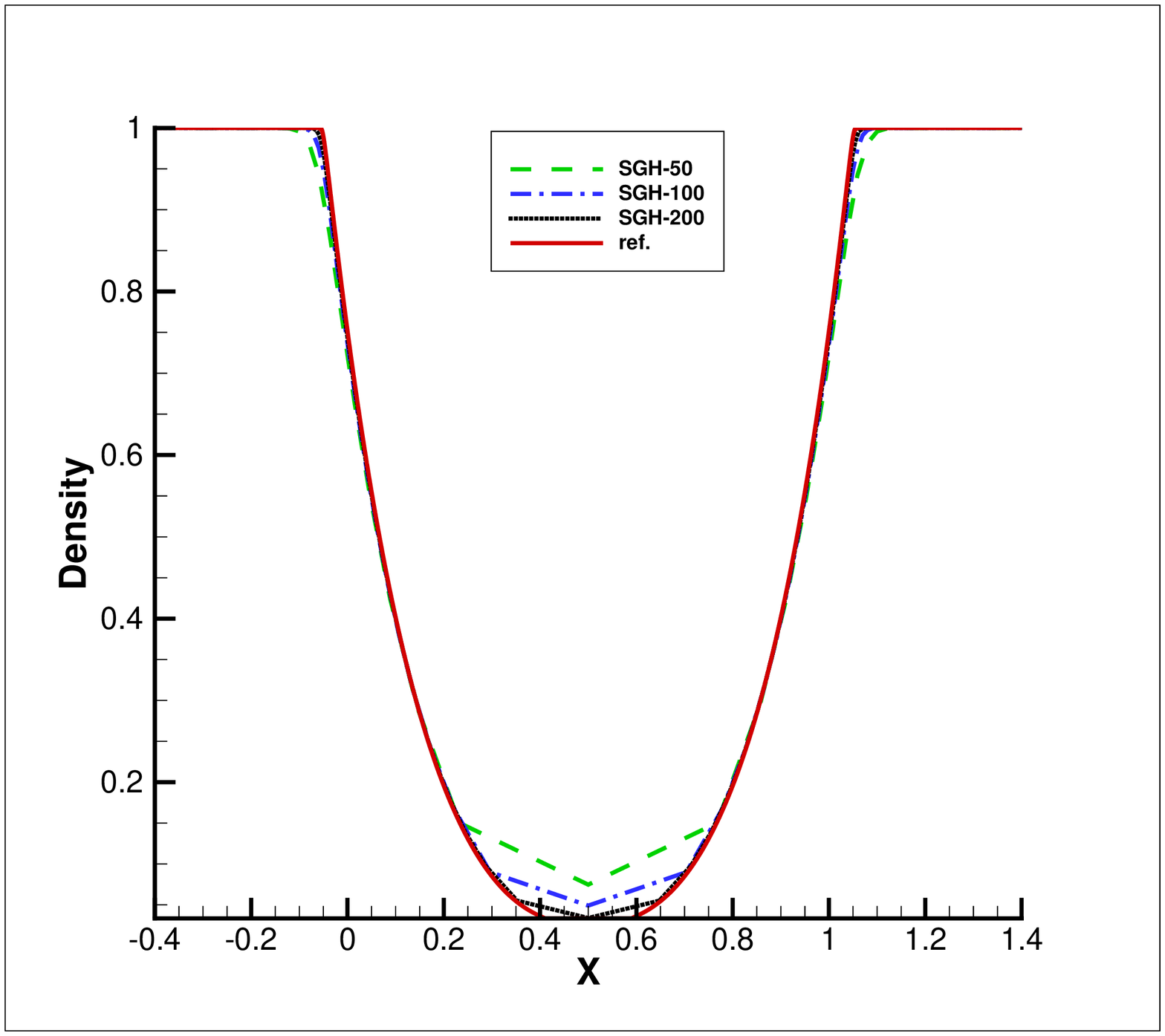}
	\includegraphics[width=2 in,keepaspectratio]{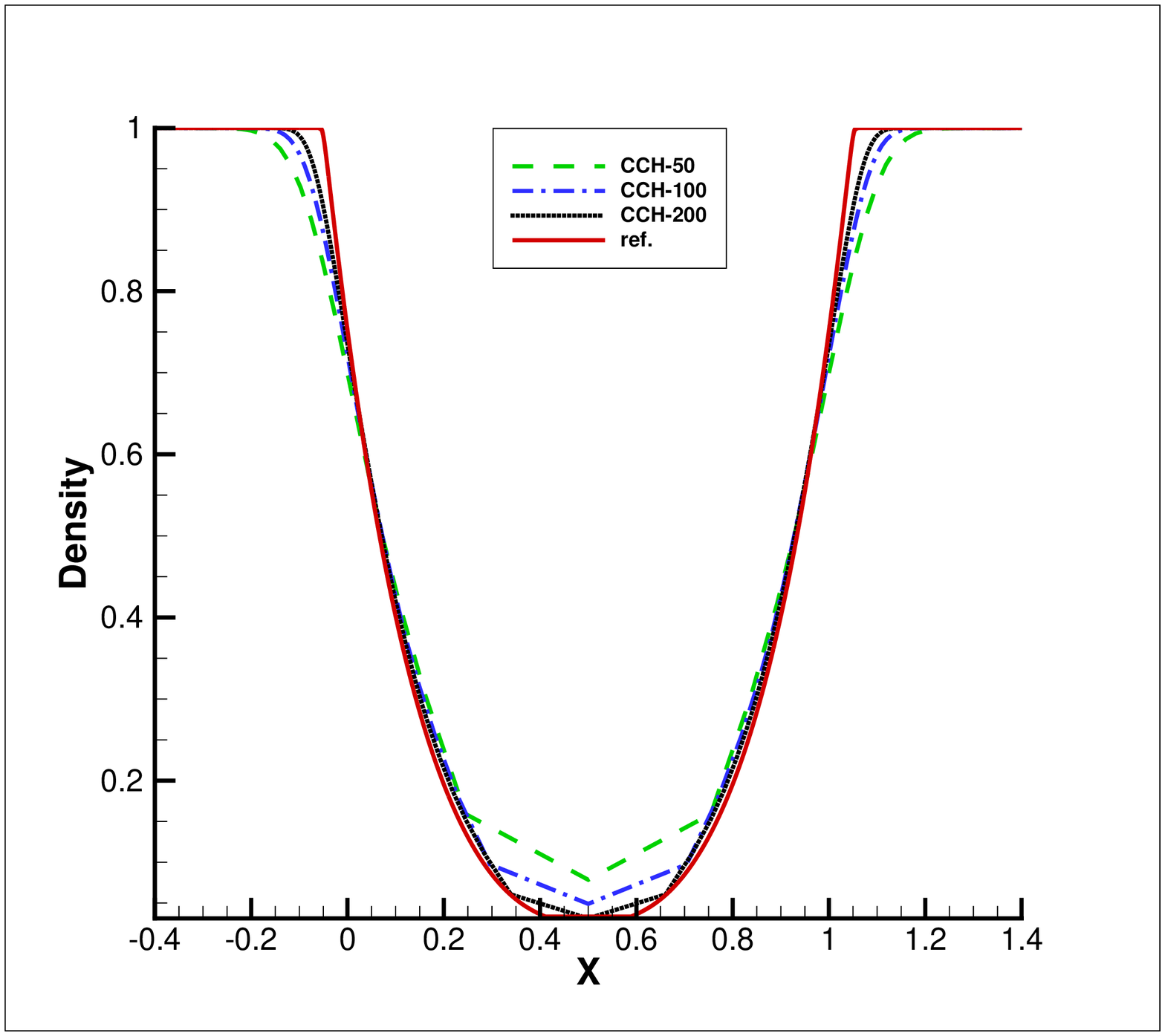}
	\caption{Double rarefaction waves: mesh convergence for density profiles with SGH (left) method and CCH (right) method, time $T=0.15$, mesh-size $N=50$ (green\ dashed\ line), $N=100$ (blue\ dashed-dotted\ line), $N=200$ (black\ dotted\ line). \label{s123-density-mj}}
\end{figure}

\begin{figure}[!h]
\centering
	\includegraphics[width=2 in,keepaspectratio]{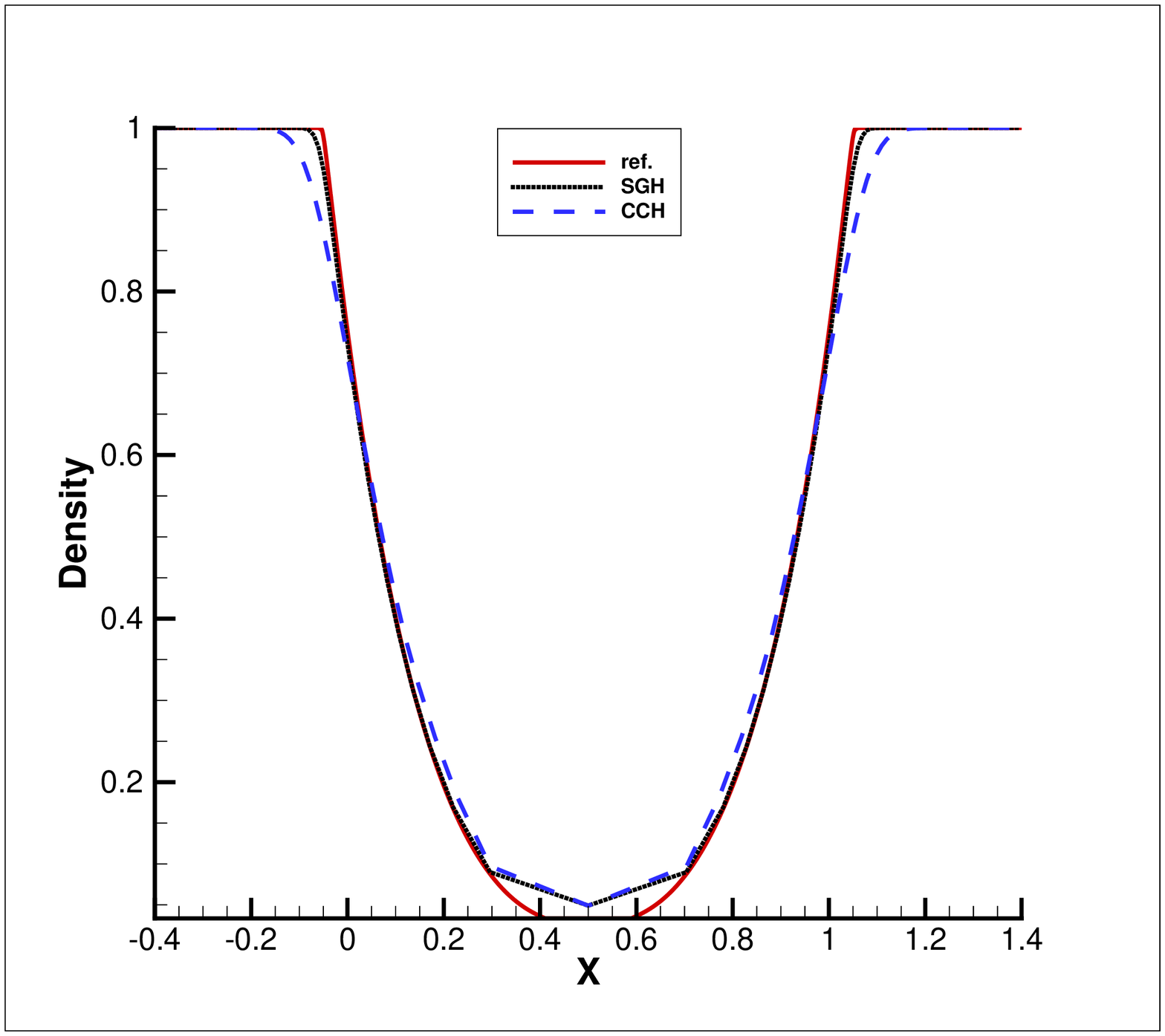}
	\includegraphics[width=2 in,keepaspectratio]{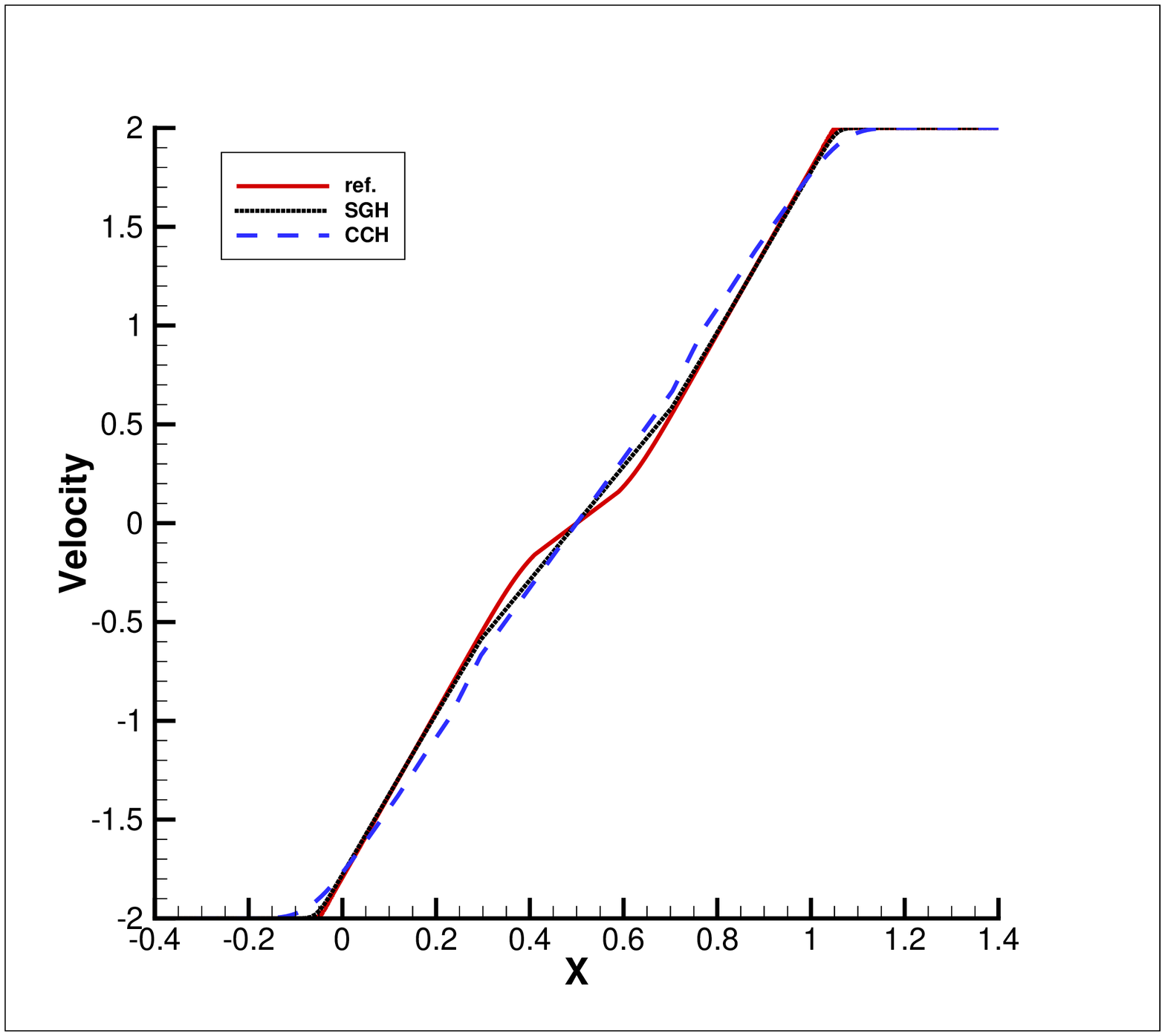}
	\caption{Double rarefaction waves: comparison between SGH (black dotted line) method and CCH (blue dashed line) method for density (left) and velocity (right) plots, time $T=0.15$, mesh-size $N=100$.\label{s123-den-vel-vs}}
\end{figure}

\begin{figure}[!h]
\centering	
	\includegraphics[width=2 in,keepaspectratio]{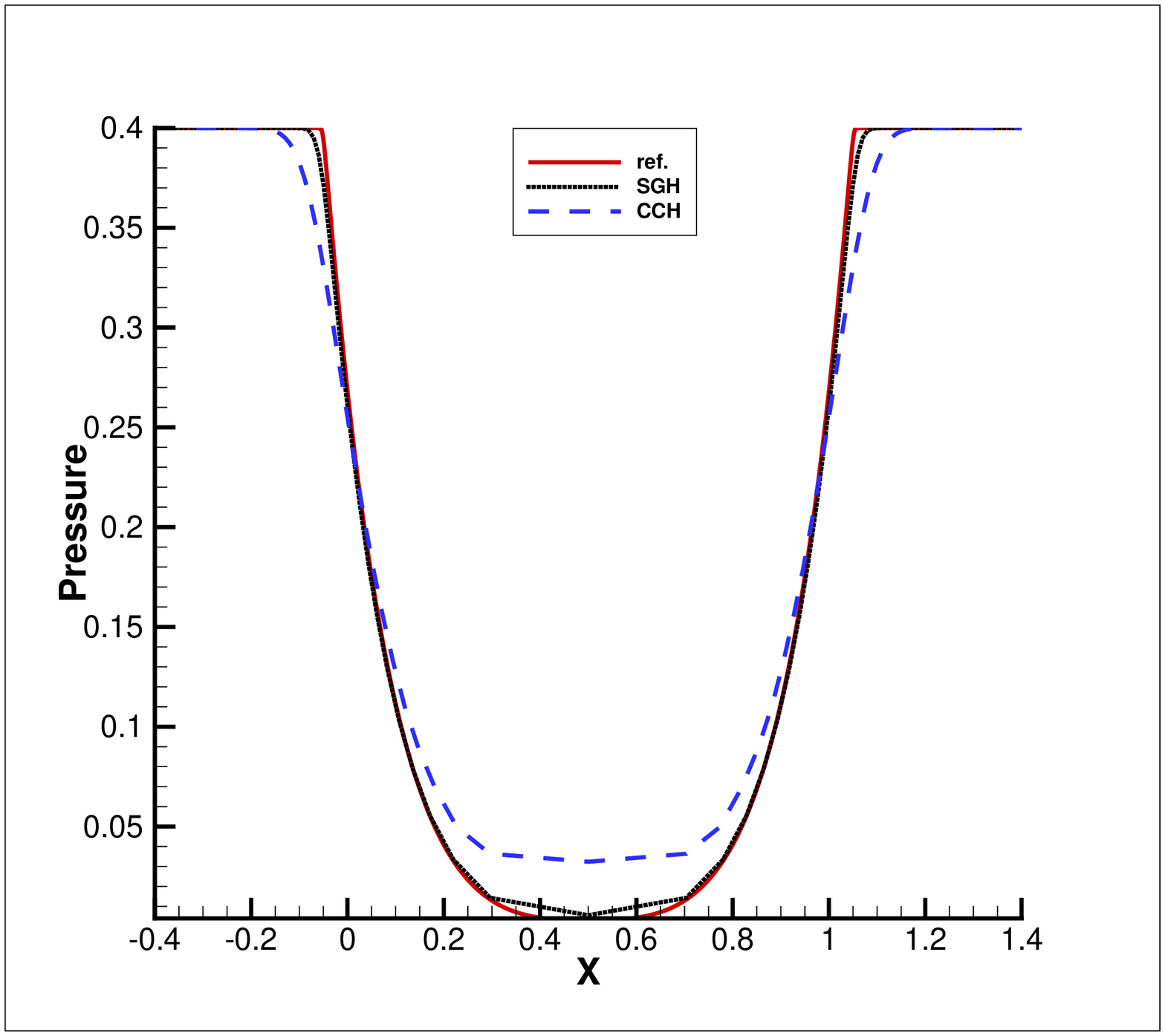}
	\includegraphics[width=2 in,keepaspectratio]{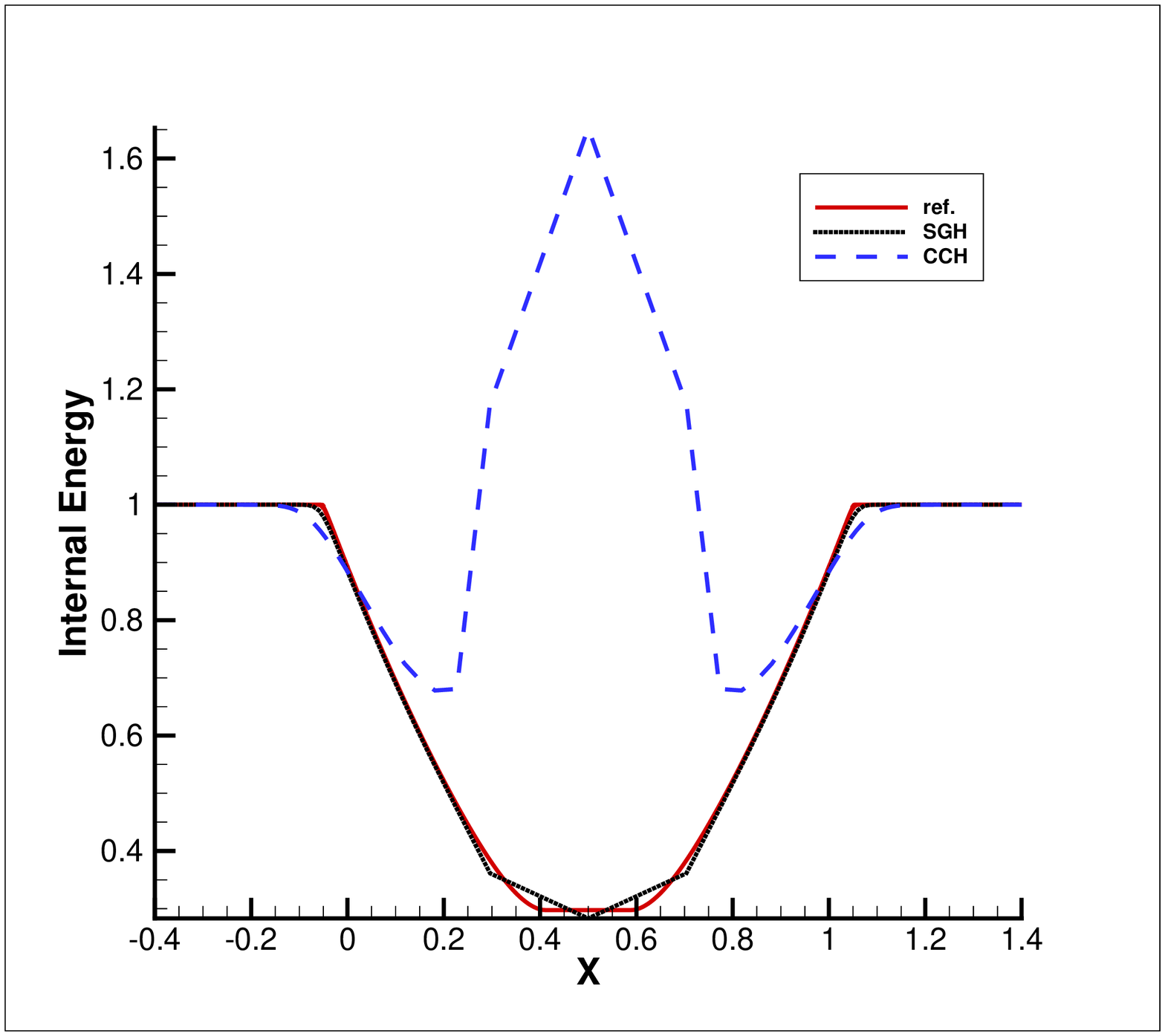}
	\caption{Double rarefaction waves: comparison between SGH (black dotted line) method and CCH (blue dashed line) method for pressure (left) and internal energy (right) plots, time $T=0.15$, mesh-size $N=100$.\label{s123-pre-ein-vs}}			
\end{figure}

Figure \ref{s123-density-mj} shows the density contours computed on different grid sizes and demonstrates that a mesh convergence study for both SGH and CCH methods leads to similar results. Similar to the previous results, it clearly shows the low-dissipative characteristic of the SGH method. In the vacuum region, there is practically no visible difference in the final state of density and velocity obtained by SCH and CCH methods, see Figure \ref{s123-den-vel-vs}. But there is some visible difference in the final state of pressure and internal energy presented in Figure \ref{s123-pre-ein-vs}.  The SGH method captures the transition between the rarefaction wave and left or right state slightly better. Although CCH method captures expansion waves at the expense of some accentuation of temperature at the center.
\subsection{Sedov blast wave problem} 
This case contains very low density with strong shocks. The computational region is $x\in[-2, 2]$ and the initial state are
$$
\left(\rho,u,E,\gamma\right)^T=\left(1,\ 0,\ {10}^{-12},\ 7/5\right)^T,
$$
everywhere except that the energy in the center cell is the constant $3.2\times10^{6}/\Delta{x}$. The inlet and outlet conditions are imposed on the left and right boundaries, respectively. The final computational time is $T =0.001$. 

\begin{figure}[!h]
\centering
	\includegraphics[width=2 in,keepaspectratio]{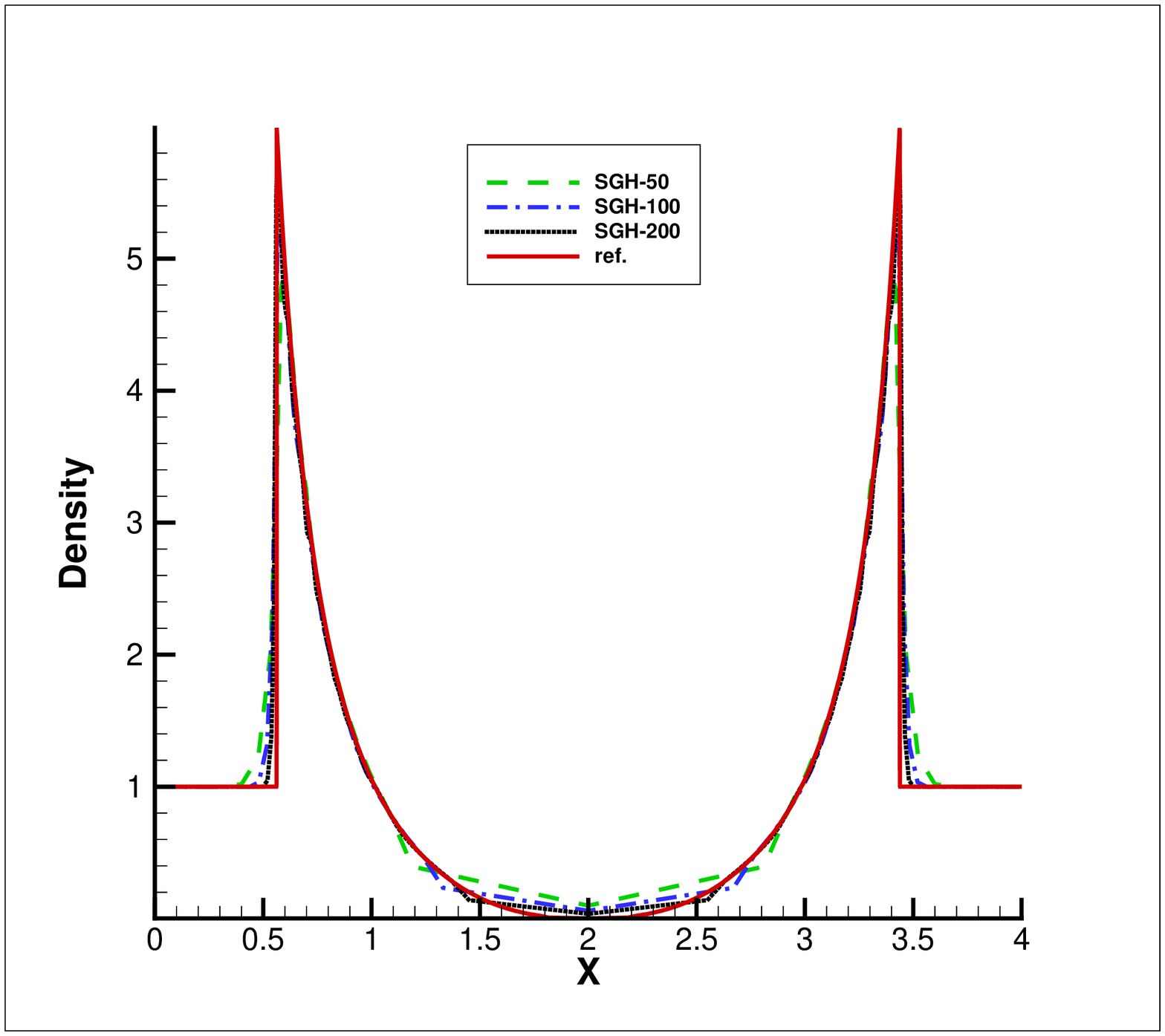}
	\includegraphics[width=2 in,keepaspectratio]{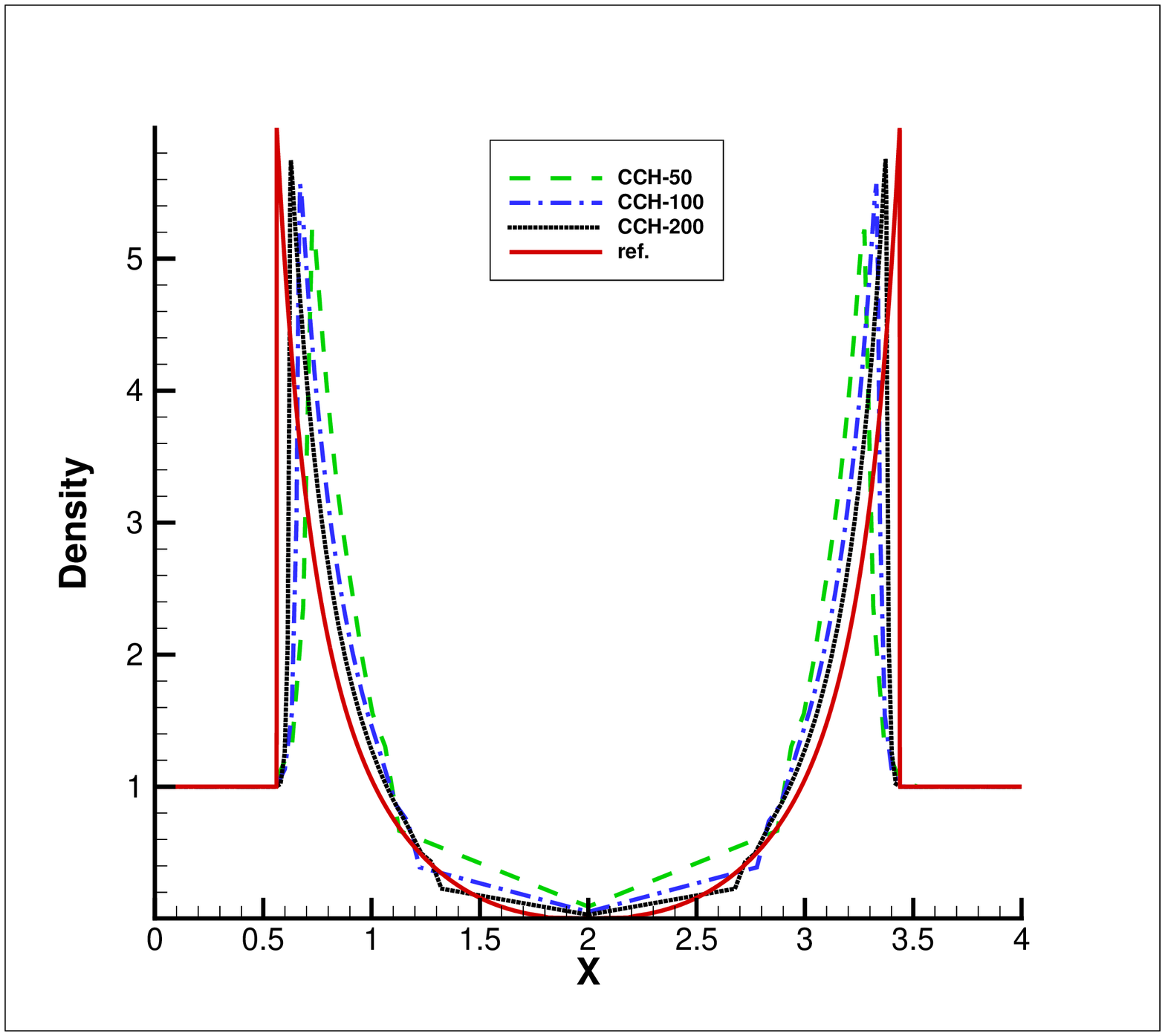}
	\caption{Sedov problem: mesh convergence for density profiles with SGH (left) method and CCH (right) method, time $T=0.001$, mesh-size $N=50$ (green\ dashed\ line), $N=100$ (blue\ dashed-dotted\ line), $N=200$ (black\ dotted\ line). \label{sedov-density-mj}}		
\end{figure}

\begin{figure}[!h]
\centering
	\includegraphics[width=2 in,keepaspectratio]{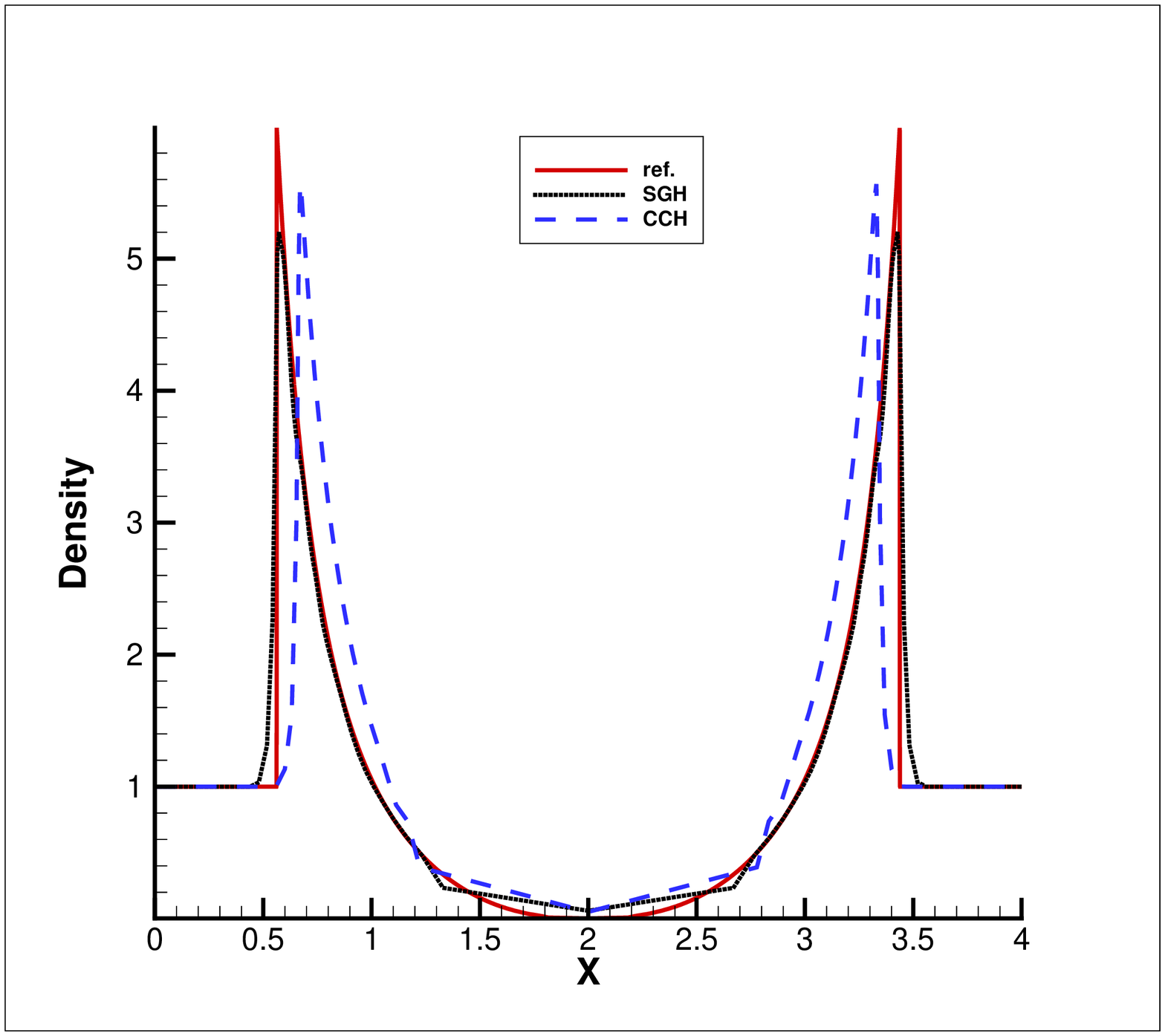}
	\includegraphics[width=2 in,keepaspectratio]{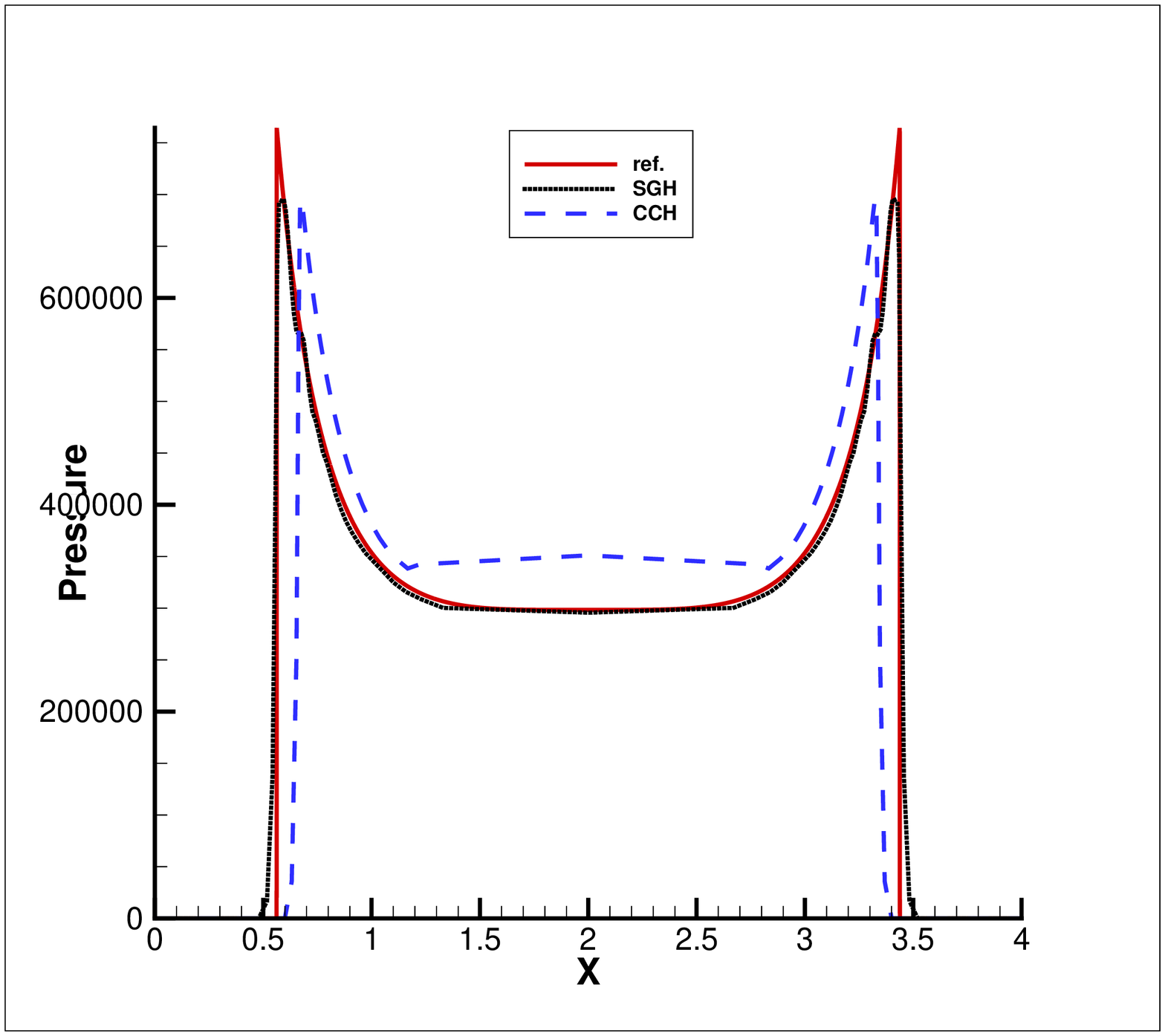}
	\caption{Sedov problem: comparison between SGH (black dotted line) method and CCH (blue dashed line) method for density (left) and pressure (right) plots, time $T=0.001$, mesh-size $N=100$.\label{sedov-den-vel-vs}}	
\end{figure}

The calculations using the acoustic solver \cite{godunov_1959} do not hold up due to the tough initial condition. However, it is observed that the CCH method with second-order approximations in (\ref{pptautau}) can solve this problem without blowing up the calculation. Figure \ref{sedov-density-mj} compares density distributions, including different mesh sizes. All of the computed solutions show little difference, and both the SGH and CCH methods exhibit excellent convergence behavior and work well for this extreme test case. The solutions of the SGH method and reference solution almost overlap each other, see Figure \ref{sedov-den-vel-vs}. The CCH method captures the discontinuities with equivalent values to the SGH method; however, it is a little behind in capturing the shock wave.
\subsection{Shock density wave interaction problem}
The simulation for the shock-entropy wave interaction problem is carried out on the domain $[-5,5]$ with the initial condition imposing natural boundary conditions
$$
\begin{array}{ll}
\left(\rho,u,P,\gamma\right)^T=\left\{
\begin{array}{ll}
\left(3.857143,\ 2.629369,\ 10.333333,\ 7/5 \right)^T  \quad &\rm{if} \quad x\leq-4,\\
\left(1+0.2\sin{5x},\ 0,\ 1,\ 7/5 \right)^T \quad &\rm{if} \quad x\geq-4.
\end{array}
\right.
\end{array}
$$
The numerical results are computed up to $T=1.8$. This test proves to be more difficult for numerical schemes since the evolution of a shock wave interacts with an entropy wave. This example assesses the capability of shock-capturing and resolving accuracy in a high-frequency local extrema region.

\begin{figure}[!h]
\centering
	\includegraphics[width=2 in,keepaspectratio]{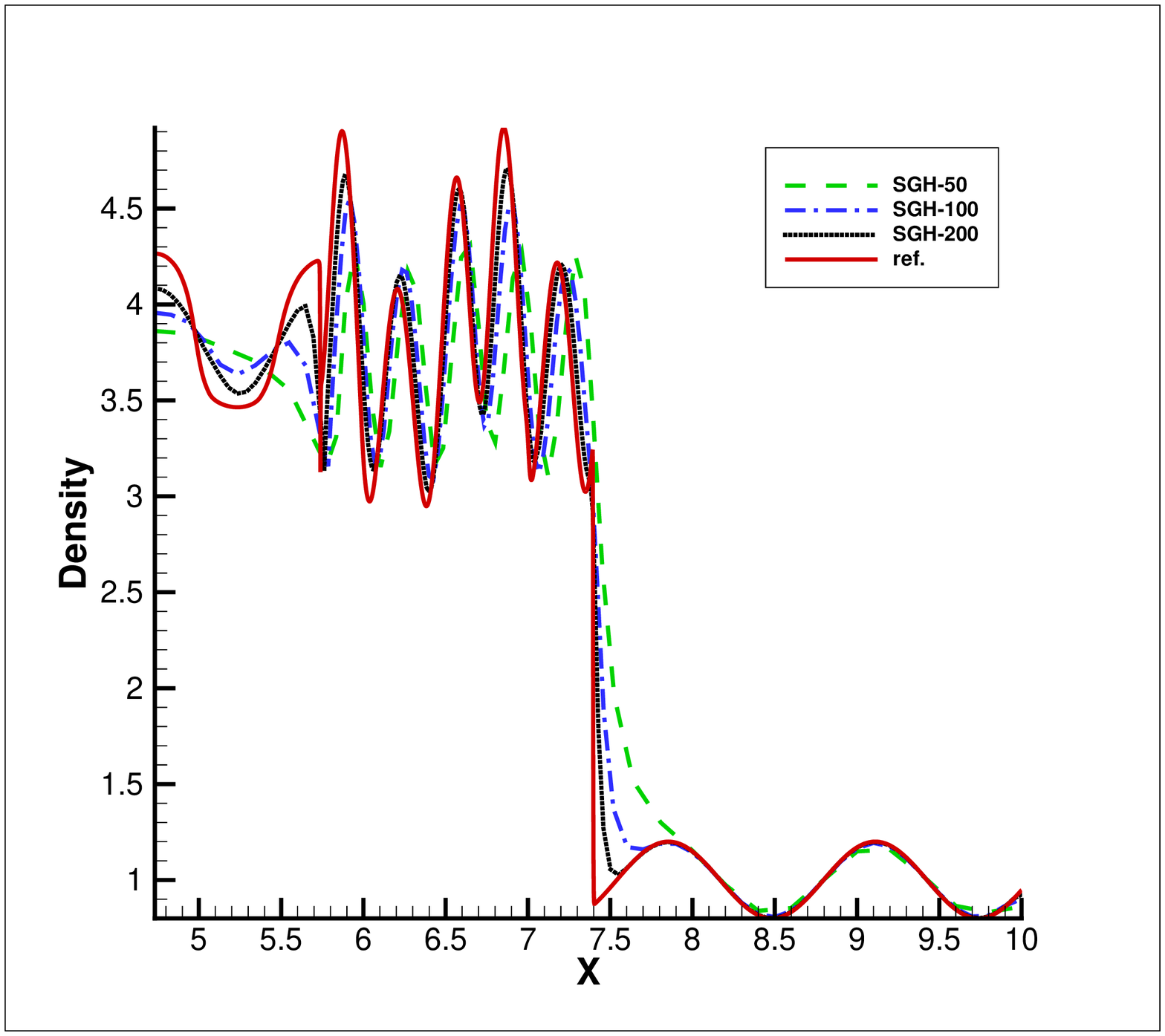}
	\includegraphics[width=2 in,keepaspectratio]{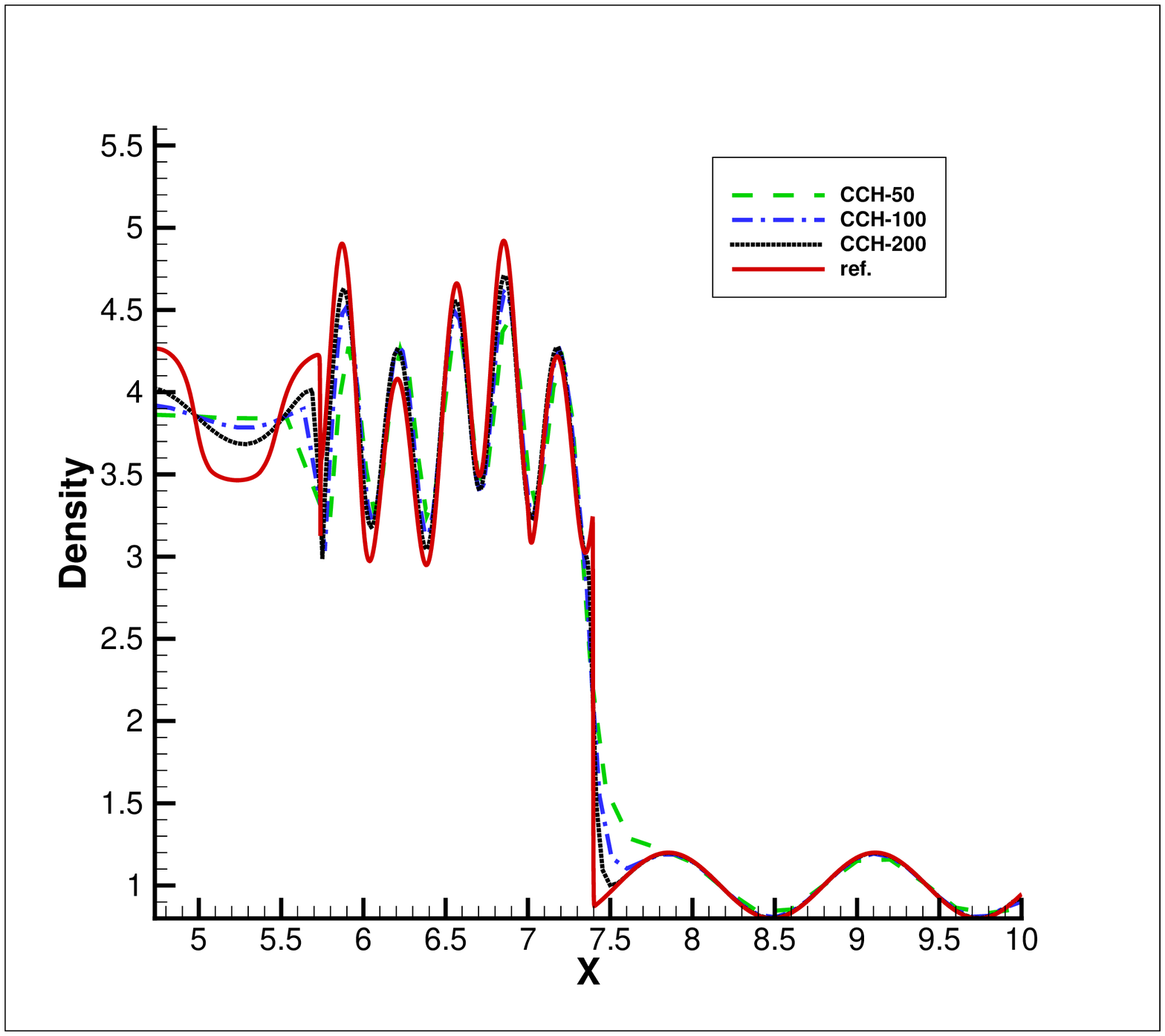}
	\caption{Shock density wave interaction problem: mesh convergence for density profiles with SGH (left) method and CCH (right) method, time $T=1.8$, mesh-size $N=50$ (green\ dashed\ line), $N=100$ (blue\ dashed-dotted\ line), $N=200$ (black\ dotted\ line).\label{shuosher-density-mj}}		
\end{figure}

\begin{figure}[!h]
\centering
	\includegraphics[width=2 in,keepaspectratio]{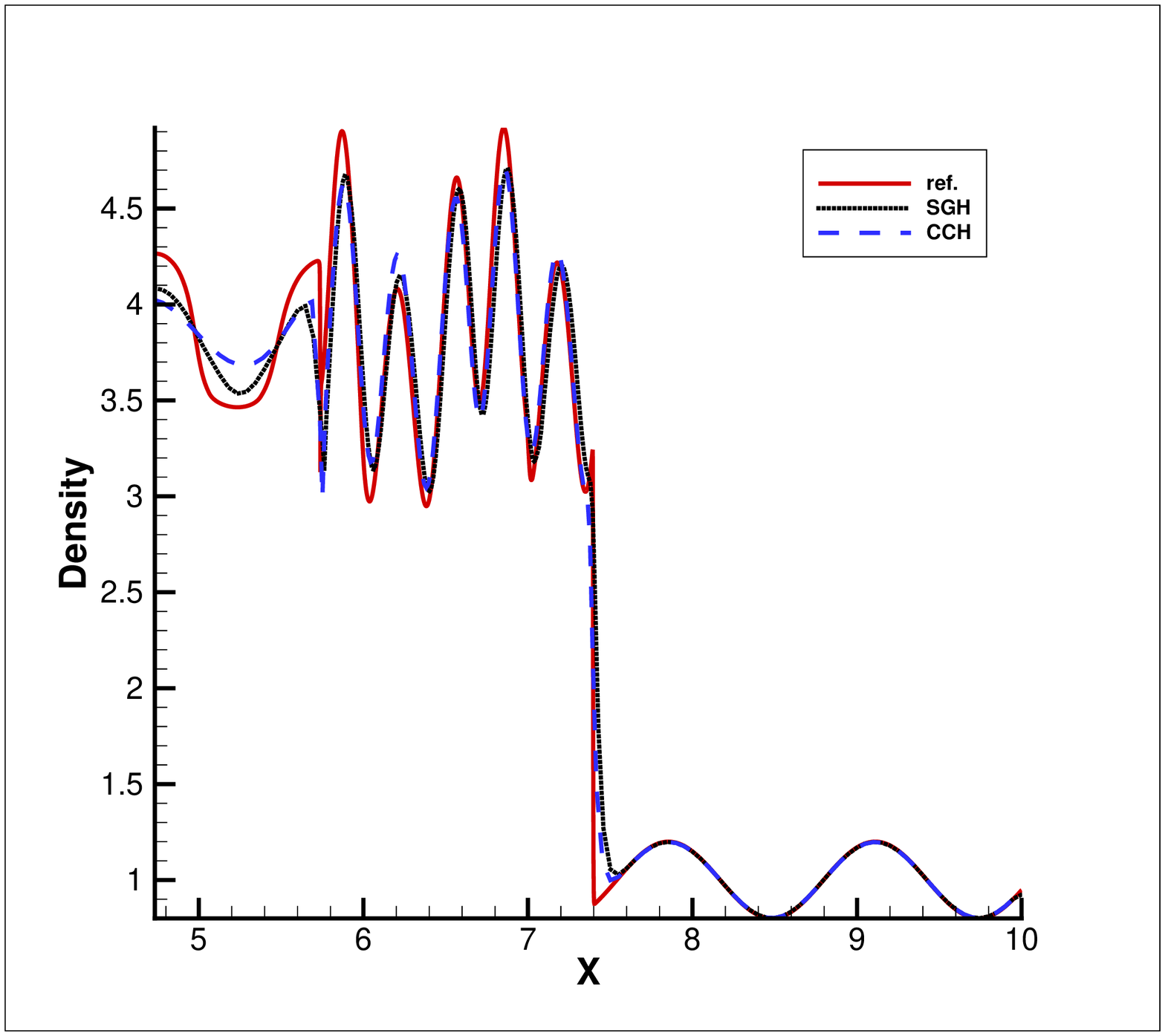}
	\includegraphics[width=2 in,keepaspectratio]{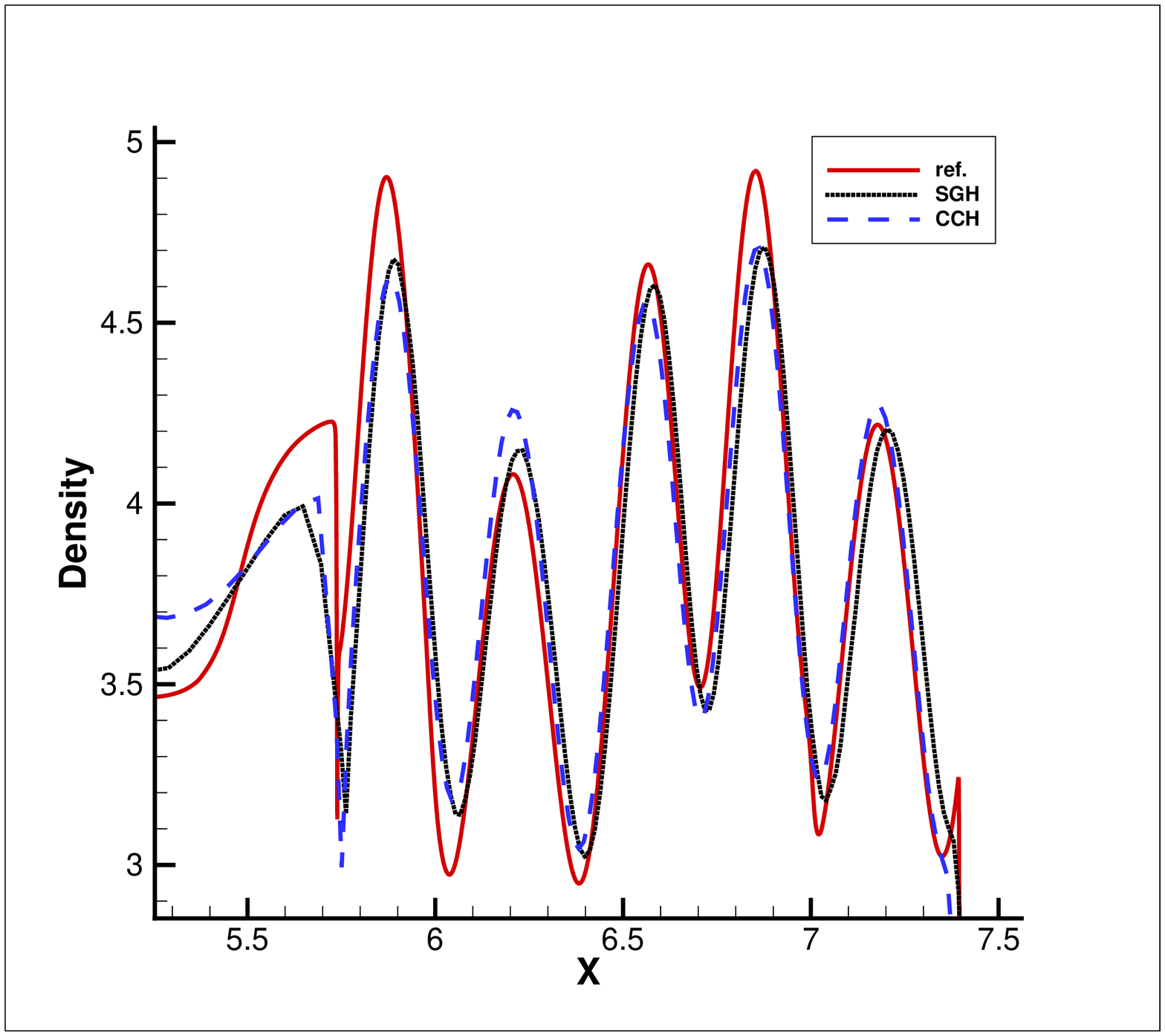}
	\caption{Shock density wave interaction problem: comparison between SGH (black dotted line) method and CCH (blue dashed line) method for density plot (left) and zoom in (right), time $T=1.8$, mesh-size $N=100$. \label{shuosher-density-vs}}	
\end{figure}

Figure \ref{shuosher-density-mj} demonstrates the solutions of density with coarse and fine grids, both results follow the overall trend of the reference solution. Figure \ref{shuosher-density-vs}  displays the density plot and its enlarged figure on the high-frequency region obtained with the SGH method; numerical computations made with the CCH method, produce similar results. Both methods have higher peaks and lower valleys and are thus closer to the reference solution. Furthermore, the SGH method does a better job of capturing the low-frequency waves more to the left of the shock. 
\subsection{LeBlanc shock tube problem}
The final problem under consideration is the LeBlanc shock tube, posed on the domain $[0,9]$, with the initial conditions as follows:
$$
\begin{array}{ll}
\left(\rho,u,P,\gamma\right)^T=\left\{
\begin{array}{ll}
\left(1,\ 0,\ {2}/{3}\times10^{-1},\ {5}/{3}\right)^T \quad &\rm{if} \quad 0<x\leq3,\\
\left(10^{-3},\ 0,\ {2}/{3}\times10^{10},\ {5}/{3}\right)^T \quad &\rm{if} \quad 3<x\leq9.
\end{array}
\right.
\end{array}
$$
A zero-gradient boundary condition is imposed on both sides. Numerical solutions are simulated at $T=6$. LeBlanc shock tube problem is the extreme case of the shock tube problem, which consists of a very strong rarefaction moving to the left and a shock wave moving to the right. These waves are separated by a huge contact wave, making this example an extraordinarily difficult numerical experiment. 

\begin{figure}[!h]
\centering
	\includegraphics[width=2 in,keepaspectratio]{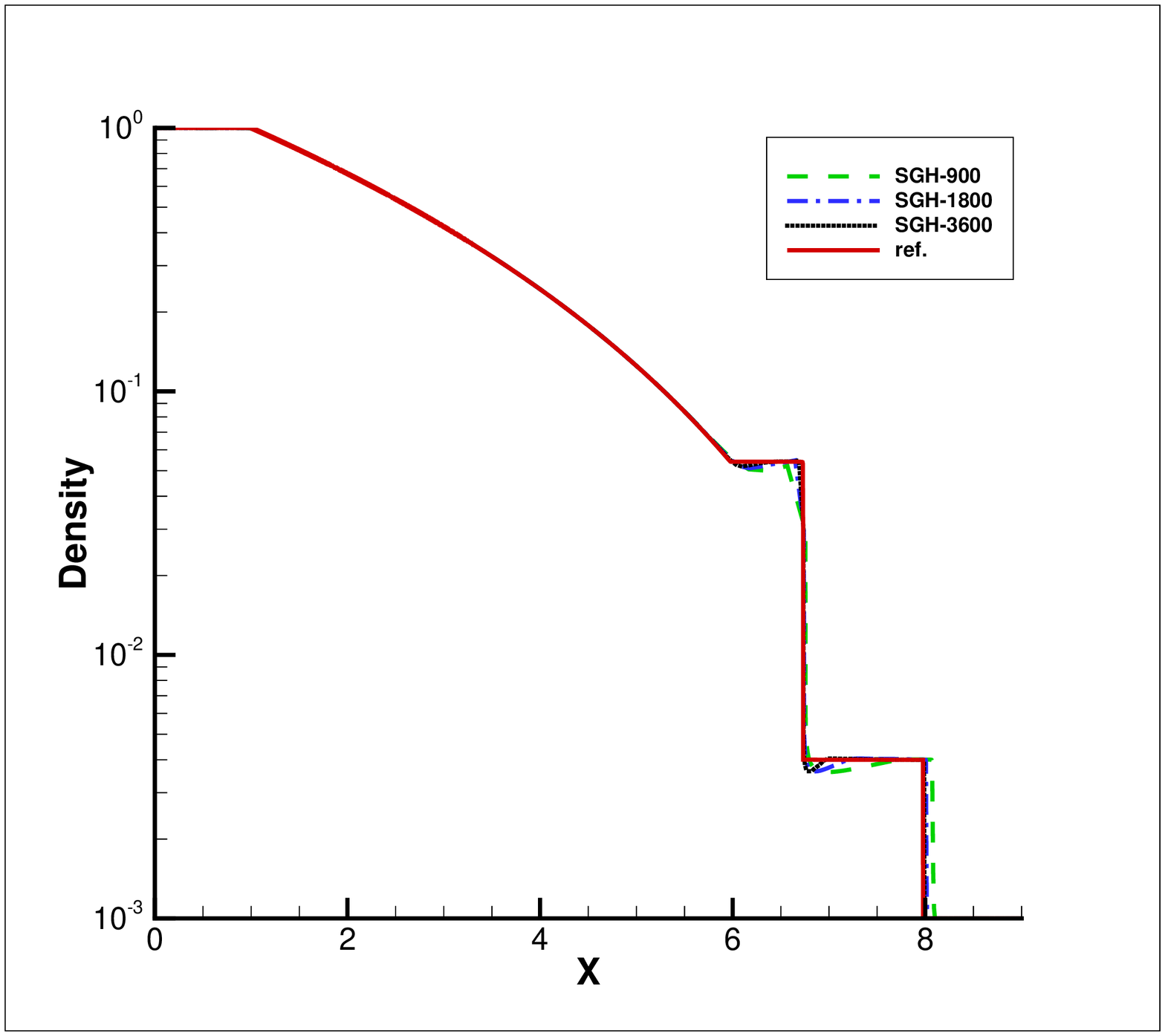}
	\includegraphics[width=2 in,keepaspectratio]{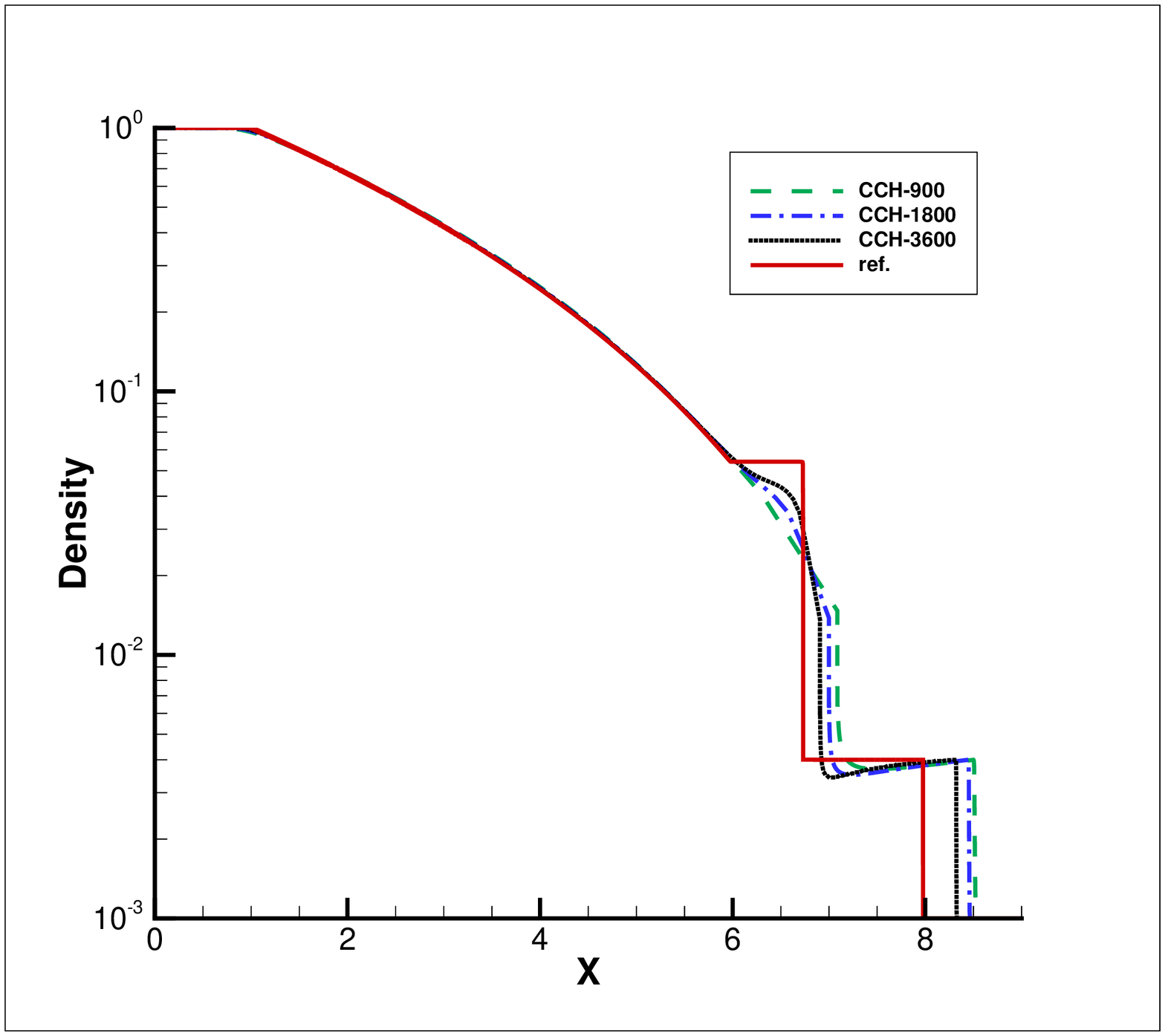}
	\caption{LeBlanc problem: mesh convergence for density profiles with SGH (left) method and CCH (right) method, time $T=6$, mesh-size $N=900$ (green\ dashed\ line), $N=1800$ (blue\ dashed-dotted\ line), $N=3600$ (black\ dotted\ line). \label{leblanc-density-mj}}	
\end{figure}

\begin{figure}[!h]
\centering
	\includegraphics[width=2 in,keepaspectratio]{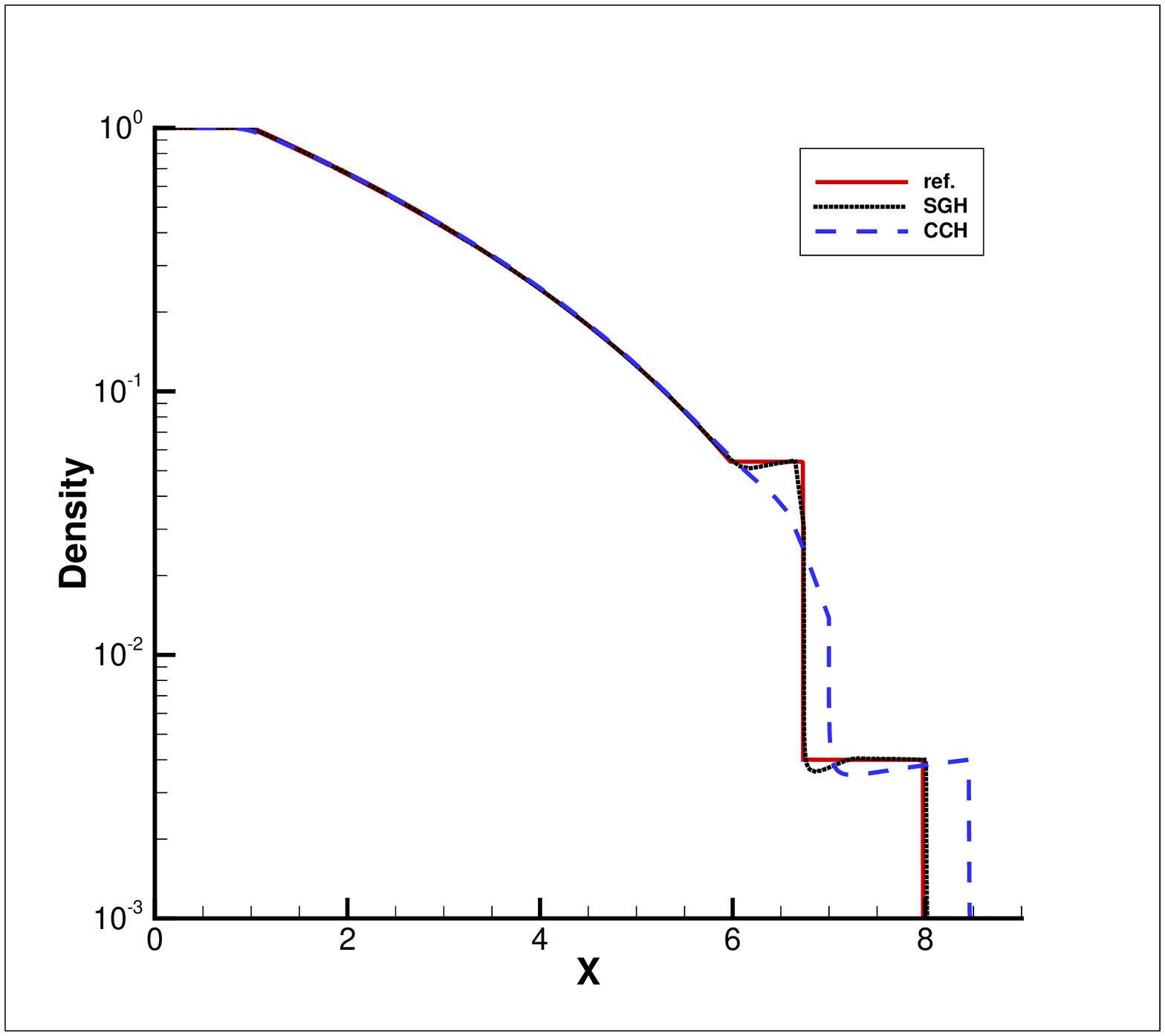}
	\includegraphics[width=2 in,keepaspectratio]{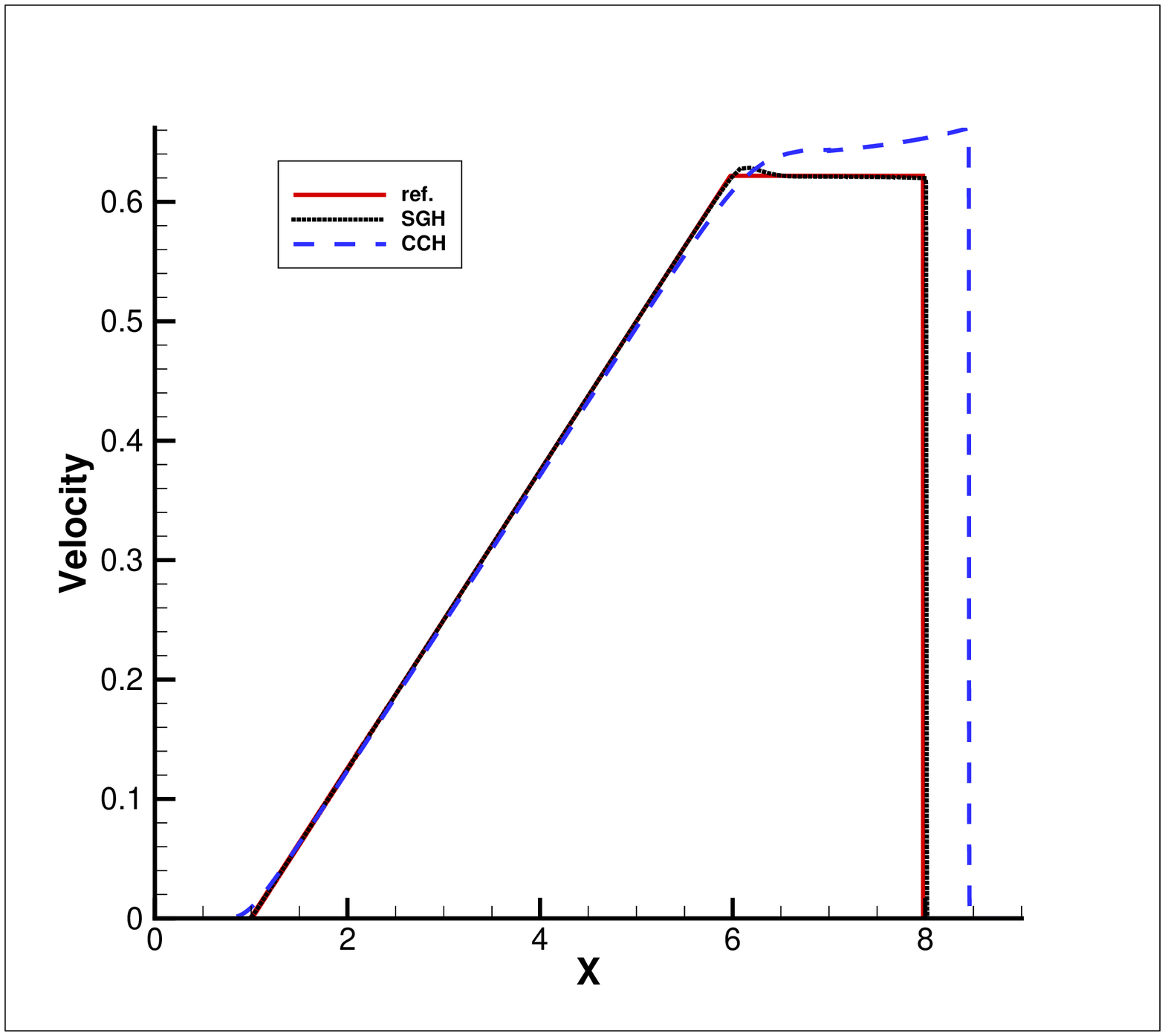}
	\caption{LeBlanc problem: comparison between SGH (black dotted line) method and CCH (blue dashed line) method for density (left) and velocity (right) plots, time $T=6$, mesh-size $N=1800$. \label{leblanc-density-vel-vs}}		
\end{figure}

\begin{figure}[!h]
\centering
	\includegraphics[width=2 in,keepaspectratio]{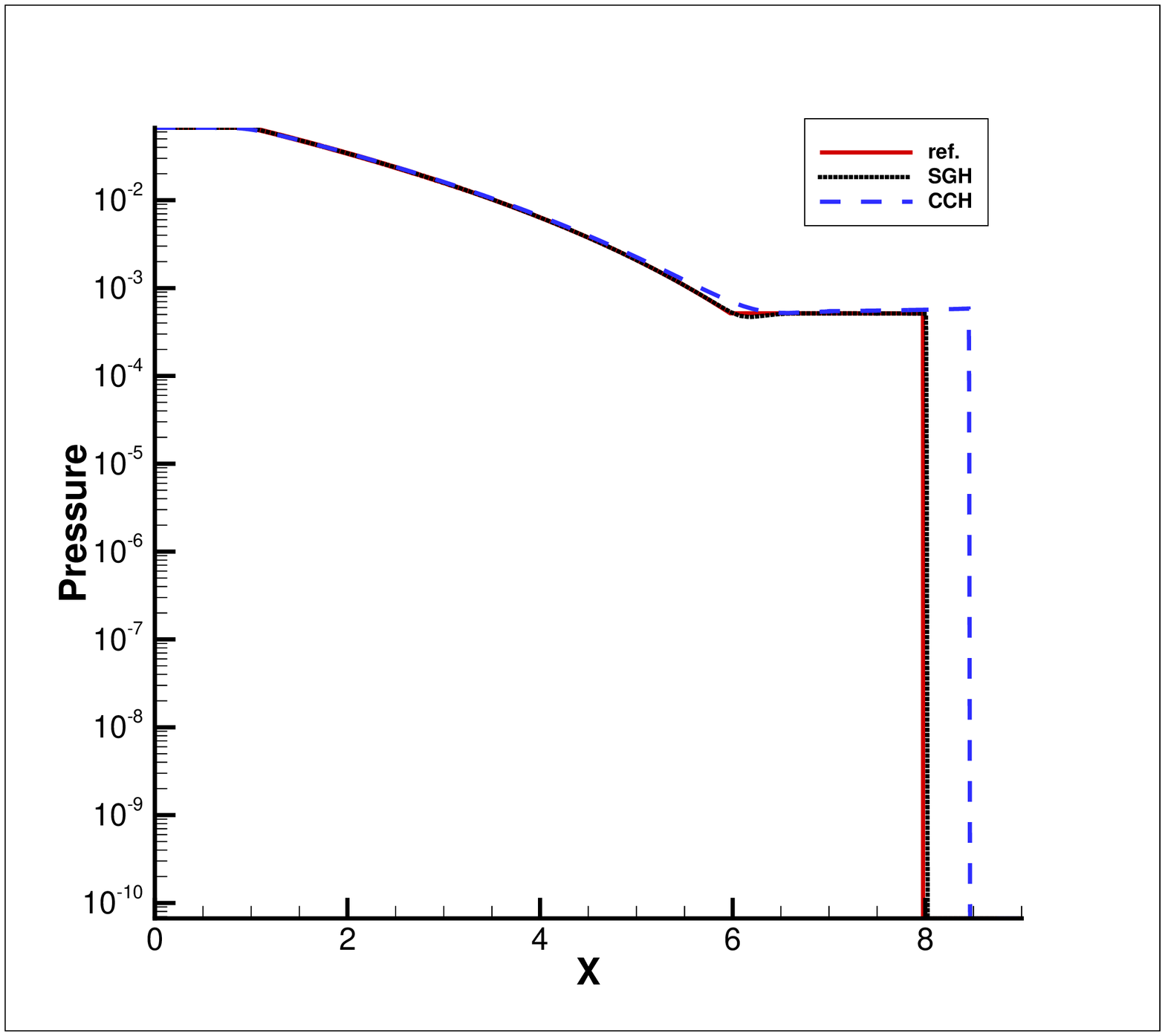}
	\includegraphics[width=2 in,keepaspectratio]{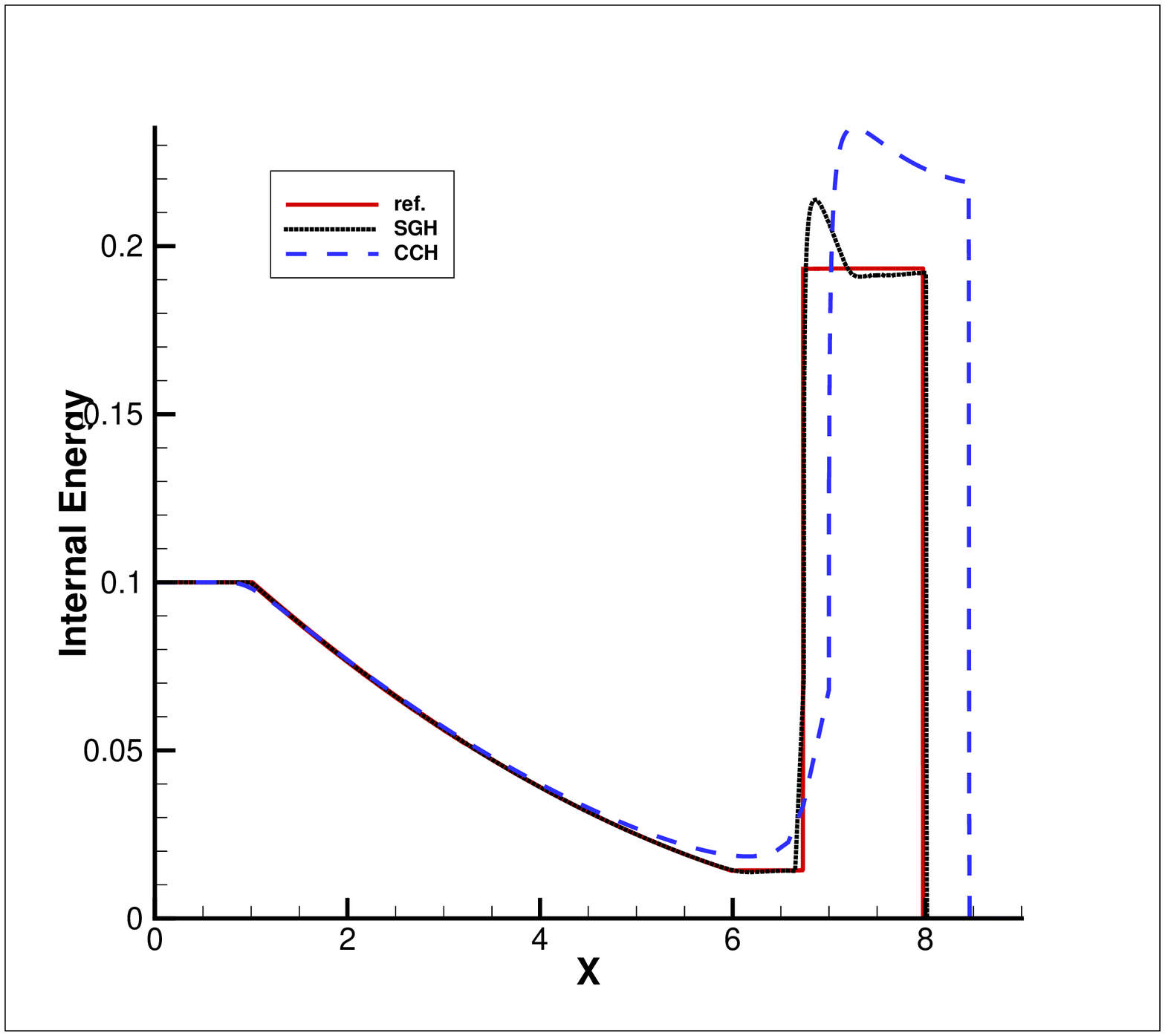}
	\caption{LeBlanc problem: comparison between SGH (black dotted line) method and CCH (blue dashed line) method for pressure (left) and internal energy (right) plots, time $T=6$, mesh-size $N=1800$. \label{leblanc-pre-ein-vs}}		
\end{figure}

Figure \ref{leblanc-density-mj} presents the results of the grid refinement study, which demonstrates the enhanced mesh convergence of both the SGH and CCH methods. Y-axis is expressed in the logarithmic scale since the density and pressure presented in this case are very small. In both cases of coarse and fine grids, the SGH method shows a better location of the shock wave than the CCH method. However, even with the finer grid, the solution of the latter method tends to overestimate the correct shock speed. The SGH method computes discontinuities that are more smeared, see Figure \ref{leblanc-density-vel-vs}. In addition, the velocity and internal energy are more accurately predicted with the SGH method in Figure \ref{leblanc-pre-ein-vs}. The SGH method tends to produce smaller overshoots in the internal energy at the contact discontinuity. These results confirm that SGH method can be an efficient numerical method to solve shock waves even in the extreme flow case.
\section{Conclusions}
In this paper, a new framework of Lagrangian schemes is developed to unify the SGH and CCH methods. The two different methods both apply this new scheme to discretize the corresponding set of equations. The scheme contains two important keys. One is the relationship between pressure and velocity, the other is Newton's second law. The main contribution of this work is that neither empirical parameters nor exact/approximate Riemann solvers are employed in this scheme. Both SGH and CCH methods using the same scheme are conservative in total mass, total momentum, and total energy and satisfy the entropy condition. Numerical results show the robustness and accuracy of the two methods. The scheme is also a new finding to build the connection between SGH and CCH methods. 

The next natural step is to generalize this strategy to multi-dimensional schemes. Note that this scheme constitutes a certain departure from the standard dimension-by-dimension approach for constructing a multidimensional scheme. The relationship between pressure and velocity can be easily extended to multi-dimension. Regarding the future research directions, a more general question should be raised in the hourglass problem. Some completed results show that the hourglass problem encountered in the SGH method rarely occurs in the CCH method. Using the same scheme, both methods may give us another way to reveal the hourglass problem from a different viewpoint.


\bibliographystyle{siamplain}
\bibliography{UfMJ}


\end{document}

%% file: ex_shared.tex
\usepackage{lipsum}
\usepackage{amsfonts}
\usepackage{graphicx}
\usepackage{epstopdf}
\usepackage{algorithmic}
\usepackage{framed}
\ifpdf
  \DeclareGraphicsExtensions{.eps,.pdf,.png,.jpg}
\else
  \DeclareGraphicsExtensions{.eps}
\fi

\newsiamremark{remark}{Remark}
\newsiamremark{hypothesis}{Hypothesis}
\crefname{hypothesis}{Hypothesis}{Hypotheses}
\newsiamthm{claim}{Claim}

\headers{Unification of Lagrangian staggered-grid hydrodynamics and cell-centered hydrodynamics in one dimension}{Xihua Xu}

\title{Unification of Lagrangian staggered-grid hydrodynamics and cell-centered hydrodynamics in one dimension}

\author{Xihua Xu \thanks{Beijing Institute of Applied Physics and Computational Mathematics, Beijing, China. (\email{xu\_xihua@iapcm.ac.cn}).}}

\usepackage{amsopn}
